\documentclass[english]{elsarticle}
\usepackage[T1]{fontenc}
\usepackage[latin9]{inputenc}
\usepackage[a4paper]{geometry}
\usepackage{euler}
\usepackage{amsmath}
\usepackage{amssymb}
\usepackage{setspace}
\usepackage{esint}
\usepackage{babel}
\begin{document}
\title{A Class of Extended Hypergeometric Functions and Its Applications}
\author{Luo Minjie \\[4pt] \footnotesize \emph{Department of Applied Mathematics, Donghua University, Shanghai 201620, PR China}\\[4pt] E-mail: \texttt{Mathwinnie@myopera.com}}
\begin{abstract}
{\normalsize Recently, there emerges different versions of beta function and hypergeometric functions containing extra parameters. Gaining enlightenment from these ideas, we will first introduce a new extension of generalized hypergeometric function and then put forward some fundamental results in the paper. Next, we will derive some properties of certain functions like extended Gauss hypergeometric functions, extended Appell's hypergeometric functions $\mathbf{F}_{1}^{\left(\kappa_{l}\right)}$, $\mathbf{F}_{2}^{\left(\kappa_{l}\right)}$, and extended Lauricella's hypergeometric function $\mathbf{F}_{D,\left(\kappa_{l}\right)}^{\left(r\right)}$,$\mathbf{F}_{A,\left(\kappa_{l}\right)}^{\left(r\right)}$, including transformation formulas, finite sum representations, Mellin-Barnes type integral representations and recurrence relations. Moreover, by using some new integral representation which will be presented in this paper, a Hardy-Hilbert type inequality involving extended Gauss hypergeometric functions will be established. }{\normalsize \par}
\end{abstract}

\maketitle

\section{{\normalsize Introduction and preliminaries.}}
In the eighteenth century, L. Euler (1707-1783) concerned himself
with the problem of interpolating between the numbers
\[
n!=\int_{0}^{\infty}e^{-t}t^{n}dt,\ \ n=0,1,2,\cdots
\]
with nonintergral values of $n$. This problem led Euler in 1729 to
the now famous gamma function, a generalization of the factorial function
that gives meaning to $x!$ when $x$ is any positive number. His
result can be extended to certain negative numbers and even to complex
numbers \cite{Special Functions for Engineers and Applied Mathematicians}. The integral representation of now widely accepted
gamma function $\Gamma\left(x\right)$ is
\begin{equation}
\Gamma\left(x\right)=\int_{0}^{\infty}t^{x-1}e^{-t}dt,\ \Re\left(x\right)>0.
\end{equation}
In 1994, by inserting a regularization factor $e^{-bt^{-1}}$, Chaudhry
and Zubair \cite{Generalized incomplete gamma functions with applications} have introduced the following extension of gamma
function
\begin{equation}
\Gamma_{p}\left(x\right)=\int_{0}^{\infty}t^{x-1}\exp\left(-t-\frac{b}{t}\right)dt,\ \Re\left(b\right)>0.
\end{equation}
For $\Re\left(b\right)>0$ this factor removes the singularity coming
from the $t=0$ limit and for $b=0$ reduces to the original gamma
function defined by (1). In 1997, Chaudhry et al. \cite{Extension of Euler's beta function} presented
the following extension of Euler's beta function
\begin{equation}
B_{b}\left(x,y\right)=\int_{0}^{1}t^{x-1}\left(1-t\right)^{y-1}\exp\left(-\frac{b}{t\left(1-t\right)}\right)dt,\ \Re\left(b\right)>0.
\end{equation}
Afterwards, Chaudhry et al. \cite{Generalized incomplete gamma functions with applications} used $B_{b}\left(x,y\right)$ to
extended Gauss hypergeometric function and Kummer confluent hypergeometric
function as follows
\begin{equation}
F_{b}\left(a,b;c;z\right)=\sum_{m=0}^{\infty}\left(a\right)_{m}\frac{B_{b}\left(b+m,c-b\right)}{B\left(b,c-b\right)}\frac{z^{m}}{m!},\ \ \left(b\geq0;\ \left|z\right|<1;\ \Re\left(c\right)>\Re\left(b\right)>0\right)
\end{equation}
\begin{equation}
\Phi_{b}\left(b;c;z\right)=\sum_{m=0}^{\infty}\frac{B_{b}\left(b+m,c-b\right)}{B\left(b,c-b\right)}\frac{z^{m}}{m!},\ \ \left(b\geq0;\ \Re\left(c\right)>\Re\left(b\right)>0\right)
\end{equation}
where $\left(a\right)_{m}$ denotes the Pochhammer symbol defined,
in terms of gamma functions, by
\[
\left(a\right)_{m}=\frac{\Gamma\left(a+m\right)}{\Gamma\left(a\right)}=\begin{cases}
1 & m=0;\ a\in\mathbb{C}\backslash\left\{ 0\right\} \\
a\left(a+1\right)\left(a+2\right)\cdots\left(a+m-1\right) & m\in\mathbb{N};\ a\in\mathbb{C}.
\end{cases}
\]
For $b=0$, functions (4) and (5) reduces to the usual hypergeometric functions.

In a similar manner, M. Ali Ozarslan and E. Ozergin \cite{Some generating relations for extended hypergeometric function via generalized fractional derivative operator} defined the extension of Appell's functions $F_{1}\left(a,b,c;d;x,y;b\right)$,
$F_{2}\left(a,b,c;d,e;x,y;b\right)$ and extended Lauricella's hypergeometric function $\mathbf{F}_{D,b}^{\left(3\right)}\left(a,b,c,d;e;x,y,z\right)$.

In 2011, E. Ozergin, M. Ali Ozarslan, and A. Altin \cite{Extension of gamma beta and hypergeometric functions} introducing the following generalizations
\begin{equation}
\Gamma_{b}^{\left(\alpha,\beta\right)}\left(x\right)=\int_{0}^{\infty}t^{x-1}{}_{1}F_{1}\left(\alpha;\beta;-t-\frac{b}{t}\right)dt\ \left(\Re\left(\alpha\right)>0,\ \Re\left(\beta\right)>0,\Re\left(p\right)>0,\Re\left(x\right)>0\right)
\end{equation}
\begin{equation}
B_{b}^{\left(\alpha,\beta\right)}\left(x,y\right)=\int_{0}^{1}t^{x-1}\left(1-t\right)^{y-1}{}_{1}F_{1}\left(\alpha;\beta;\frac{-b}{t\left(1-t\right)}\right)dt.
\end{equation}
\[
\left(\Re\left(\alpha\right)>0,\ \Re\left(\beta\right)>0,\Re\left(p\right)>0,\Re\left(x\right)>0,\Re\left(y\right)>0\right)
\]

In this paper, our extension mainly based on using the following generalization of gamma and beta functions.
\newtheorem{Definition}{Definition}[section]
\begin{Definition}
\cite[p.243]{A Class of Extended Fractional Derivative Operators and Associated Generating Relations Involving Hypergeometric Functions}: Let a function $\Theta\left(\left\{ \kappa_{l}\right\} _{l\in\mathbb{N}_{0}};z\right)$
be analytic within the disk $\left|z\right|<R$ $\left(0<R<\infty\right)$
and let its Taylor-Maclaurin coeficients be explicitly deonted by
the sequence $\left\{ \kappa_{l}\right\} _{l\in\mathbb{N}_{0}}$.
Suppose also that the function $\Theta\left(\left\{ \kappa_{l}\right\} _{l\in\mathbb{N}_{0}};z\right)$
can be continued analytically in the right half-plane $\Re\left(z\right)>0$
with the asymptotic property given as follows:
\begin{equation}
\Theta\left(\kappa_{l};z\right)\equiv\Theta\left(\left\{ \kappa_{l}\right\} _{l\in\mathbb{N}_{0}};z\right)=\begin{cases}
\sum_{l=0}^{\infty}\kappa_{l}\frac{z^{l}}{l!} & \left(\left|z\right|<R;\ 0<R<\infty;\kappa_{0}=1\right)\\
M_{0}z^{\omega}\exp\left(z\right)\left[1+O\left(\frac{1}{z}\right)\right] & \left(\Re\left(z\right)\rightarrow\infty;M_{0}>0;\omega\in\mathbb{C}\right)
\end{cases}
\end{equation}
for some suitable constants $M_{0}$ and $\omega$ depending essentially on the sequence $\left\{ \kappa_{l}\right\} _{l\in\mathbb{N}_{0}}$. We can define extended Gamma function $\Gamma_{b}^{\left(\kappa_{l} \right)}\left(z\right)$ and the extended Beta function $\mathcal{B}^{\left(\kappa_{l}\right)}\left(\alpha,\beta;b\right)$ by
\begin{equation}
\Gamma_{b}^{\left(\kappa_{l}\right)}\left(z\right)=\int_{0}^{\infty}t^{z-1}\Theta\left(\left\{ \kappa_{l}\right\} ;-t-\frac{b}{t}\right)dt\ \left(\Re\left(z\right)>0;\ \Re\left(b\right)\geq0\right)
\end{equation}
and
\begin{equation}
\mathcal{B}_{b}^{\left(\kappa_{l}\right)}\left(\alpha,\beta;b\right)=\int_{0}^{1}t^{\alpha-1}\left(1-t\right)^{\beta-1}\Theta\left(\left\{ \kappa_{l}\right\} ;-\frac{b}{t\left(1-t\right)}\right)dt
\end{equation}
\[
\left(\min\left\{ \Re\left(\alpha\right),\Re\left(\beta\right)\right\} >0;\ \Re\left(b\right)\geq0\right).
\]
By introducing one additional parameter $d$ with $\Re\left(d\right)\geq0$,
we have
\begin{equation}
\mathcal{B}_{b,d}^{\left(\kappa_{l}\right)}\left(\alpha,\beta\right)=\int_{0}^{1}t^{\alpha-1}\left(1-t\right)^{\beta-1}\Theta\left(\left\{ \kappa_{l}\right\} ;-\frac{b}{t}-\frac{d}{1-t}\right)dt
\end{equation}
\[
\left(\min\left\{ \Re\left(\alpha\right),\Re\left(\beta\right)\right\} >0;\ \min\left\{ \Re\left(b\right),\Re\left(d\right)\right\} \geq0\right)
\]
\end{Definition}
\newenvironment{Remark}{{\noindent\bfseries Remark}\par}{\hfill  $\diamondsuit$\par}
\begin{Remark}
\begin{enumerate}[(I)]
\item For given sequence $\left\{ \kappa_{l}\right\} _{l\in\mathbb{N}_{0}}$our
definitions would reduce to known or new extensions of the Gamma and
Beta functions. Recall that the asymptotic behavior of Kummer's confluent
hypergeometric function at infinity has the form
\[
\Phi\left(a;c;z\right)={}_{1}F_{1}\left(a;c;z\right)=\frac{\Gamma\left(c\right)}{\Gamma\left(a\right)}e^{z}z^{a-c}\left[1+O\left(\frac{1}{z}\right)\right],\ \Re\left(z\right)\rightarrow\infty.
\]
 It is clear that $\Phi\left(a;c;z\right)$ is a special case of $\Theta\left(\kappa_{l};z\right)$ with sequence $\left\{ \frac{\left(a\right)_{l}}{\left(c\right)_{l}}\right\} _{l\in\mathbb{N}_{0}}$. And our Gamma function $\Gamma_{b}^{\left( \kappa_{l}\right)}$ and Beta function $\mathcal{B}^{\left(\kappa_{l}\right)}$
are reduced to (6) and (7) which has been studied in \cite{Extension of gamma beta and hypergeometric functions}. We can also replace $_{1}F_{1}\left(a;c;z\right)$ with $\exp\left(z\right)$ to get more basic extension of Gamma and Beta function. (see, for details, \cite{Extended hypergeometric and confluent hypergeometric functions})
\item Function $\mathcal{B}_{b,d}^{\left(\kappa_{l}\right)}\left(\alpha,\beta\right)$
is a further extension of Beta function $\mathcal{B}_{b}^{\left(\kappa_{l}\right)}\left(\alpha,\beta\right)$
by introducing complex parameter $d$ $\left(\Re\left(d\right)\geq0\right)$.
This improvement may be essential since $-\frac{b}{t}-\frac{d}{1-t}=-\left(\frac{b+\left(d-b\right)t}{t\left(1-t\right)}\right)$.
Note that if we set $d=b$ we will get function $\mathcal{B}_{b}^{\left(\kappa_{l}\right)}\left(\alpha,\beta\right)$.
In general, form $-\frac{b}{t}-\frac{d}{1-t}$ may give us more flexibility
for making transformations.
\end{enumerate}
\end{Remark}

By using $\mathcal{B}_{b,d}^{\left(\kappa_{l}\right)}\left(\alpha,\beta\right)$ we can get the corresponding further extension of Gauss hypergeometric function as follows.
\begin{Definition}
\cite{A Class of Extended Fractional Derivative Operators and Associated Generating Relations Involving Hypergeometric Functions} The Extended Gauss hypergeometric function $_{2}F_{1}^{\left(\kappa_{l}\right)}$ is defined by
\begin{equation}
_{2}F_{1}^{\left(\kappa_{l}\right)}\left[\begin{array}{c}
\alpha_{1},\alpha_{2}\\
\beta_{1}
\end{array};z;b,d\right]=\sum_{n=0}^{\infty}\left(\alpha_{1}\right)_{n}\frac{\mathcal{B}_{b,d}^{\left(\kappa_{l}\right)}\left(\alpha_{2}+n,\beta_{1}-\alpha_{2}\right)}{B\left(\alpha_{2},\beta_{1}-\alpha_{2}\right)}\frac{z^{n}}{n!}.
\end{equation}
\[
\left(\left|z\right|<1;\Re\left(c\right)>\Re\left(b\right)>0;\min\left\{ \Re\left(b\right),\Re\left(d\right)\right\} \geq0\right)
\]
\end{Definition}

If $\Theta\left(\kappa_{l};z\right)=\exp\left(z\right)$, we write
\[
_{2}F_{1}\left[\begin{array}{c}
\alpha_{1},\alpha_{2}\\
\beta_{1}
\end{array};z;b,d\right]=\sum_{n=0}^{\infty}\left(\alpha_{1}\right)_{n}\frac{\mathcal{B}_{b,d}\left(\alpha_{2}+n,\beta_{1}-\alpha_{2}\right)}{B\left(\alpha_{2},\beta_{1}-\alpha_{2}\right)}\frac{z^{n}}{n!}.
\]

Its integral representation can be easily derived.
\newtheorem{Theorem}{Theorem}[section]
\begin{Theorem}
For the extended Gauss hypergeometric function $_{2}F_{1}^{\left(\kappa_{l}\right)}\left(\alpha_{1},\alpha_{2};\beta_{1};z;b,d\right)$,
we have the following integral represetation:
\begin{equation}
_{2}F_{1}^{\left(\kappa_{l}\right)}\left[\begin{array}{c}
\alpha_{1},\alpha_{2}\\
\beta_{1}
\end{array};z;b,d\right]=\frac{1}{B\left(\alpha_{2},\beta_{1}-\alpha_{2}\right)}\int_{0}^{1}t^{\alpha_{2}-1}\left(1-t\right)^{\beta_{1}-\alpha_{2}-1}\left(1-zt\right)^{-\alpha_{1}}\Theta\left(\kappa_{l};-\frac{b}{t}-\frac{d}{1-t}\right)dt,
\end{equation}
\[
\Re\left(\beta_{1}\right)>\Re\left(\alpha_{2}\right)>0;\  \Re\left(b\right),\Re\left(d\right)>0; b=d=0, \left|\arg\left(1-z\right)<\pi\right|
\]
\end{Theorem}

To prove this integral representation, we first expand $\left(1-zt\right)^{-\alpha_{1}}$ as
its Taylor series
\[
\left(1-zt\right)^{-\alpha_{1}}=\sum_{n=0}^{\infty}\left(\alpha_{1}\right)_{n}\frac{\left(zt\right)^{n}}{n!}.
\]
It is clear that this series is absolutely and uniformly convergent. Thus, we can interchange the order of summation and than integrate out $t$ to get the result.\\[5pt]

As outlined above, we have introduced the general idea of the development path of the subject, the extended Beta function and its application in generalizing some hypergeometric series. As we all know, general hypergeometric functions have important applications in mathematical, physical and chemical areas. Thus, quite naturally, we should not neglect the question that whether we could extend the definition of generalized hypergeometric function by using extended beta function $\mathcal{B}_{b,d}^{\left(\kappa_{l}\right)}\left(\alpha,\beta\right)$. Or whether we could prove that general hypergeometric functions are specific cases of the functions we get in the paper. The significance of such questions lies in the fact that a large number of practical problems require more complex hypergeometric functions of complicated coefficients. Fortunately, the answers are almost certainly "yes", inspired by the way how extended Appell's hypergeometric functions $\mathbf{F}_{2}^{\left(\kappa_{l}\right)}\left(\alpha,\beta_{1},\beta_{2};\gamma_{1},\gamma_{2};x,y;b,d\right)$ are constructed in \cite{A Class of Extended Fractional Derivative Operators and Associated Generating Relations Involving Hypergeometric Functions}, \cite{Some generating relations for extended hypergeometric function via generalized fractional derivative operator}. In fact, the method can also be applied to ratiocinate a new kind of multivariate hypergeometric function, i.e., $\mathbf{F}_{A,\left(\kappa_{l}\right)}^{\left(r\right)}\left(\alpha,\beta_{1},\cdots,\beta_{r};\gamma_{1},\cdots\gamma_{r};x_{1},\cdots,x_{r};b,d\right)$, with its large numbers of formulas, which will be analyzed in section 3.

Let us review some basic conclusions concerning with general hypergeometric functions. The generalized hypergeometric function with $p$ numerator and $q$
denominator parameters is defined by \cite[p. 27]{Theory and Applications of Fractional Differential Equations}
\begin{align}
_{p}F_{q}\left(\alpha_{1},\cdots,\alpha_{p};\beta_{1},\cdots,\beta_{q};z\right)={}_{p}F_{q}\left[\begin{array}{ccc}
\alpha_{1} & \cdots & \alpha_{p}\\
\beta_{1} & \cdots & \beta_{q}
\end{array};z\right] & =\sum_{k=0}^{\infty}\frac{\left(\alpha_{1}\right)_{k}\cdots\left(\alpha_{p}\right)_{k}}{\left(\beta_{1}\right)_{k}\cdots\left(\beta_{q}\right)_{k}}\frac{z^{k}}{k!},\\
 & =\frac{\Gamma\left(\beta_{1}\right)\cdots\Gamma\left(\beta_{q}\right)}{\Gamma\left(\alpha_{1}\right)\cdots\Gamma\left(\alpha_{p}\right)}\sum_{k=0}^{\infty}\frac{\Gamma\left(\alpha_{1}+k\right)\cdots\Gamma\left(\alpha_{p}+k\right)}{\Gamma\left(\beta_{1}+k\right)\cdots\Gamma\left(\beta_{q}+k\right)}\frac{z^{k}}{k!}
\end{align}
\[
\left(\alpha_{l},\beta_{j}\in\mathbb{C},\beta_{j}\neq0,-1,-2,\cdots,l=1,\cdots,p;j=1,\cdots,q\right)
\]
which is absolutely convergent for all values of $z\in\mathbb{C}$ if $p\leq q$. When $p=q+1$, the series is absolutely convergent for $\left|z\right|<1$ and for $\left|z\right|=1$ when $\Re\left(\sum_{j=1}^{q}\beta_{j}-\sum_{l=1}^{p}\alpha_{l}\right)>0$,
while it is conditionally convergent for $\left|z\right|=1$ $\left(z\neq1\right)$
if $-1<\Re\left(\sum_{j=1}^{q}\beta_{j}-\sum_{l=1}^{p}\alpha_{l}\right)\leq0$. More detailed information may be found in \cite{Special Functions for Engineers and Applied Mathematicians, Special Functions}and \cite{SPECIAL FUNCTIONS}.

In what follows, we consider the extended generalized hypergeometric
function $_{p}F_{q}^{\left(\kappa_{l}\right)}\left(\alpha_{1},\cdots\alpha_{p};\beta_{1}\cdots\beta_{q};z;b,d\right)$
defined by
\begin{equation}
_{p}F_{q}^{\left(\kappa_{l}\right)}\left[\begin{array}{ccc}
\alpha_{1} & \cdots & \alpha_{p}\\
\beta_{1} & \cdots & \beta_{q}
\end{array};z;b,d\right]=
\begin{cases}
\displaystyle\sum_{m=0}^{\infty}\left(\alpha_{1}\right)_{m}\prod_{j=1}^{q}\frac{\mathcal{B}_{b,d}^{\left(\kappa_{l}\right)}\left(\alpha_{j+1}+m,\beta_{j}-\alpha_{j+1}\right)}{B\left(\alpha_{j+1},\beta_{j}-\alpha_{j+1}\right)}\frac{z^{m}}{m!}, & \left|z\right|<1;\ p=q+1\\[15pt]
\displaystyle\sum_{m=0}^{\infty}\prod_{j=1}^{q}\frac{\mathcal{B}_{b,d}^{\left(\kappa_{l}\right)}\left(\alpha_{j}+m,\beta_{j}-\alpha_{j}\right)}{B\left(\alpha_{j},\beta_{j}-\alpha_{j}\right)}\frac{z^{m}}{m!}, & p=q\\[15pt]
\displaystyle\sum_{m=0}^{\infty}\frac{1}{\left(\beta_{1}\right)_{m}}\cdots\frac{1}{\left(\beta_{r}\right)_{m}}\prod_{j=1}^{p}\frac{\mathcal{B}_{b,d}^{\left(\kappa_{l}\right)}\left(\alpha_{j}+m,\beta_{r+j}-\alpha_{j}\right)}{B\left(\alpha_{j},\beta_{r+j}-\alpha_{j}\right)}\frac{z^{m}}{m!}, & r=q-p,\ p<q.
\end{cases}
\end{equation}
If we set $b=d=0$, then
\[
_{p}F_{q}^{\left(\kappa_{l}\right)}\left(\alpha_{1},\cdots\alpha_{p};\beta_{1}\cdots\beta_{q};z;b,d\right)={}_{p}F_{q}\left(\alpha_{1},\cdots,\alpha_{p};\beta_{1},\cdots,\beta_{q};z\right).
\]
If $p=2,\ q=1$ then $_{p}F_{q}^{\left(\kappa_{l}\right)}$ reduces
to extended Gauss hypergeometric function (12).\\[10pt]

The paper is organized as follows.

In section 2, we will prove Euler type integral representation and expound some properties of new extended generalized hypergeometric function $_{p}F_{q}^{\left(\kappa_{l}\right)}$. By using \emph{Ramanujan's Master Theorem}, we can get the Mellin-Barnes type integral representation for extended generalized hypergeometric function $_{p}F_{q}^{\left(\kappa_{l}\right)}$. Then, as a special case of $_{p}F_{q}^{\left(\kappa_{l}\right)}$, more properties of extended Gauss hypergeometric functions such as differentiation formulas and recurrence formulas will be discussed in detail. In the last part of this section, we give an improved definition, \textbf{definition 2.2}. By doing this, we are not only able to prove a summation theorem deserving special attention, but also to establish the relations between modified functions and fractional operators.

In section 3, we will give more integral identities and reduction formulas for functions such as extended Appell's hypergeometric functions $\mathbf{F}_{1}^{\left(\kappa_{l}\right)}\left(\alpha,\beta,\beta';\gamma;x,y;b,d\right)$,
$\mathbf{F}_{2}^{\left(\kappa_{l}\right)}\left(\alpha,\beta,\beta';\gamma,\gamma';x,y;b,d\right)$ and Laucrella's hypergeometric functions $\mathbf{F}_{D,\left(\kappa_{l}\right)}^{\left(r\right)}$, defined in \cite{A Class of Extended Fractional Derivative Operators and Associated Generating Relations Involving Hypergeometric Functions}. In addition, we give a definition of a new type of multivariate hypergeometric function, i.e.,
\begin{multline}
\mathbf{F}_{A,\left(\kappa_{l}\right)}^{\left(r\right)}\left(\alpha,\beta_{1},\cdots,\beta_{r};\gamma_{1},\cdots\gamma_{r};x_{1},\cdots,x_{r};b,d\right)=\\ \sum_{m_{1},\cdots, m_{r}=0}^{\infty}\left(\alpha\right)_{m_{1}+\cdots+m_{r}}\prod_{j=1}^{r}\frac{\mathcal{B}_{b,d}^{\left(\kappa_{l}\right)}\left(\beta_{j}+m_{j},\gamma_{j}-\beta_{j}\right)}{B\left(\beta_{j},\gamma_{j}-\beta_{j}\right)}\frac{x_{1}^{m_{1}}\cdots x_{r}^{m_{r}}}{m_{1}!\cdots m_{r}!}
\end{multline}
\[
\left(\left|x_{1}\right|+\cdots+\left|x_{r}\right|<1;\ \min\left\{ \Re\left(b\right),\Re\left(d\right)\right\} \geq0\right)
\]
and draw conclusions about it on the basis of the definition. Viewed as a a generalization of \emph{Ramanujan's Master Theorem} and with its wide application in evaluating multidimensional definite integrals and the Feynman integral, the \emph{Method of Bracket} is also used to prove the Mellin-Barnes integral representations of these multivariate hypergeometric functions.

In section 4, we reestablish a Hilbert-Hardy inequality by using extended Gauss hypergeometric functions. It is remarkable that the extended Gauss hypergeometric functions can also be used to generalize other inequalities.\\[10pt]

\section{{\normalsize Extended Generalized Hypergeometric Functions}}
We first give a formal definition and some remarks for extended generalized hypergeometric functions.
\begin{Definition}
For suitably constrained (real or complex) parameters
$\alpha_{j},\ j=1,\cdots p$; $\beta_{i},\ i=1,\cdots,q$, we define
the extended generalized hypergeometric functions by
\begin{multline}
_{p}F_{q}^{\left(\kappa_{l}\right)}\left[\begin{array}{ccc}
\alpha_{1} & \cdots & \alpha_{p}\\
\beta_{1} & \cdots & \beta_{q}
\end{array};z;b,d\right]=\\
\begin{cases}
\displaystyle\sum_{m=0}^{\infty}\left(\alpha_{1}\right)_{m}\prod_{j=1}^{q}\frac{\mathcal{B}_{b,d}^{\left(\kappa_{l}\right)}\left(\alpha_{j+1}+m,\beta_{j}-\alpha_{j+1}\right)}{B\left(\alpha_{j+1},\beta_{j}-\alpha_{j+1}\right)}\frac{z^{m}}{m!}, & \left|z\right|<1;\ p=q+1; \Re(\beta_j)>\Re(\alpha_{j+1})>0\\[15pt]
\displaystyle\sum_{m=0}^{\infty}\prod_{j=1}^{q}\frac{\mathcal{B}_{b,d}^{\left(\kappa_{l}\right)}\left(\alpha_{j}+m,\beta_{j}-\alpha_{j}\right)}{B\left(\alpha_{j},\beta_{j}-\alpha_{j}\right)}\frac{z^{m}}{m!}, & z\in\mathbb{C}; p=q; \Re(\beta_j)>\Re(\alpha_j)>0\\[15pt]
\displaystyle\sum_{m=0}^{\infty}\frac{1}{\left(\beta_{1}\right)_{m}}\cdots\frac{1}{\left(\beta_{r}\right)_{m}}\prod_{j=1}^{p}\frac{\mathcal{B}_{b,d}^{\left(\kappa_{l}\right)}\left(\alpha_{j}+m,\beta_{r+j}-\alpha_{j}\right)}{B\left(\alpha_{j},\beta_{r+j}-\alpha_{j}\right)}\frac{z^{m}}{m!}, & z\in\mathbb{C}; r=q-p,\ p<q; \Re(\beta_{r+j}>\Re(\alpha_j)>0)
\end{cases}
\end{multline}
\end{Definition}
\begin{Remark}
\begin{enumerate}[(I)]
\item If the numerator parameter $\alpha_{1}$ in the first series
of (18) is zero or a negative integer, the series terminates, i.e.,
\[
_{q+1}F_{q}^{\left(\kappa_{l}\right)}\left[\begin{array}{ccc}
-n & \cdots & \alpha_{q+1}\\
\beta_{1} & \cdots & \beta_{q}
\end{array};z;b,d\right]=\sum_{m=0}^{n}\left(-n\right)_{m}\prod_{j=1}^{q}\frac{\mathcal{B}_{b,d}^{\left(\kappa_{l}\right)}\left(\alpha_{j+1}+m,\beta_{j}-\alpha_{j+1}\right)}{B\left(\alpha_{j+1},\beta_{j}-\alpha_{j+1}\right)}\frac{z^{m}}{m!}.
\]
And we should make sure that no denominator parameter $\beta_{1},\cdots,\beta_{r}$ is allowed to be zero or a negative integer in the third expression
of (18).
\item The convergence of series $_{p}F_{q}^{\left(\kappa_{l}\right)}$
with different parameter $p$ and $q$ can be obtained by using Comparison
Test and the convergence of known corresponding generalized hypergeometric
functions $_{p}F_{q}$ . Their own convergence region is therefore
clear.
\item Note that for usual generalized hypergeometric function
$_{p}F_{q}\left(\alpha_{1},\cdots,\alpha_{p};\beta_{1},\cdots,\beta_{q};z\right)$
, its numerator-parameters $\alpha_{1},\cdots\alpha_{p}$ and denominator-parameters
$\beta_{1},\cdots,\beta_{q}$ can be arbitrarily rearranged without
any effects on its value. But under the setting of our definition,
we need to point out that the expression shown on the right-hand side
of (18) is closely associated with the order of the parameters given
in the square bracket on the the left-hand side of (18).
\item Although the general definition seems complicated, it will
be very convenient in practice.
\end{enumerate}
\end{Remark}

We will first discuss in detail some fundamental properties of function (18) and
then in the final part of this section a modified definition will be presented.\\

\subsection{{\normalsize \textbf{Some basic results about the extended generalized hypergeometric
function}}}
What we most concerned about is the Euler type integral representation of our new function (18). The following theorem aims to demonstrate that the form of  the Euler type integral representation of $_{p}F_{q}^{\left(\kappa_{l}\right)}$
is very similar to that of the Euler type integral representation of $_{p}F_{q}$.
\begin{Theorem}
If $p\leq q+1$, $\Re\left(\beta_{q}\right)>\Re\left(\alpha_{p}\right)>0$
and $\Re\left(b\right)>0,\ \Re\left(d\right)>0$, we have
\begin{align}
_{p}F_{q}^{\left(\kappa_{l}\right)}\left[\begin{array}{ccc}
\alpha_{1} & \cdots & \alpha_{p}\\
\beta_{1} & \cdots & \beta_{q}
\end{array};z;b,d\right] & =\frac{\Gamma\left(\beta_{q}\right)}{\Gamma\left(\alpha_{p}\right)\Gamma\left(\beta_{q}-\alpha_{p}\right)}\int_{0}^{1}t^{\alpha_{p}-1}\left(1-t\right)^{\beta_{q}-a_{p}-1}\nonumber \\
 & \ \ \times{}_{p-1}F_{q-1}^{\left(\kappa_{l}\right)}\left[\begin{array}{ccc}
\alpha_{1} & \cdots & \alpha_{p-1}\\
\beta_{1} & \cdots & \beta_{q-1}
\end{array};zt;b,d\right]\Theta\left(\kappa_{l};-\frac{b}{t}-\frac{d}{1-t}\right)dt.
\end{align}
Eq (19) also holds when $b=d=0$, proided that $\left|\arg\left(1-z\right)<\pi\right|$.
\end{Theorem}
\newenvironment{Proof}{{\noindent\bfseries Proof}\par}{\hfill $\square$\par}
\begin{Proof}
We need to verify that formula (19) holds for three different expressions of $_{p}F_{q}^{\left(\kappa_{l}\right)}\left(\alpha_{1},\cdots,\alpha_{p};\beta_{1},\cdots,\beta_{q};z;b,d\right)$ given in (18), respectively. Consider the case $p=q+1$. Note that
\[
\frac{\mathcal{B}_{b,d}^{\left(\kappa_{l}\right)}\left(\alpha_{q+1}+m,\beta_{q}-\alpha_{q+1}\right)}{B\left(\alpha_{q+1},\beta_{q}-\alpha_{q+1}\right)}=\frac{\Gamma\left(\beta_{q}\right)}{\Gamma\left(\alpha_{q+1}\right)\Gamma\left(\beta_{q}-\alpha_{q+1}\right)}\int_{0}^{1}t^{\alpha_{q+1}+m-1}\left(1-t\right)^{\beta_{q}-\alpha_{q+1}-1}\Theta\left(\kappa_{l};-\frac{b}{t}-\frac{d}{1-t}\right)dt.
\]
Substituting this in expression (18), we get
\begin{multline}
_{q+1}F_{q}^{\left(\kappa_{l}\right)}\left[\begin{array}{ccc}
\alpha_{1} & \cdots & \alpha_{q+1}\\
\beta_{1} & \cdots & \beta_{q}
\end{array};z;b,d\right] =\sum_{m=0}^{\infty}\left(\alpha_{1}\right)_{m}\prod_{j=1}^{q-1}\frac{\mathcal{B}_{b,d}^{\left(\kappa_{l}\right)}\left(\alpha_{j+1}+m,\beta_{j}-\alpha_{j+1}\right)}{B\left(\alpha_{j+1},\beta_{j}-\alpha_{j+1}\right)}
\frac{\Gamma\left(\beta_{q}\right)}{\Gamma\left(\alpha_{q+1}\right)\Gamma\left(\beta_{q}-\alpha_{q+1}\right)}\\\int_{0}^{1}t^{\alpha_{q+1}+m-1}\left(1-t\right)^{\beta_{q}-\alpha_{q+1}-1}\Theta\left(\kappa_{l};-\frac{b}{t}-\frac{d}{1-t}\right)dt\cdot\frac{z^{m}}{m!}\\
=\frac{\Gamma\left(\beta_{q}\right)}{\Gamma\left(\alpha_{q+1}\right)\Gamma\left(\beta_{q}-\alpha_{q+1}\right)}\int_{0}^{1}t^{\alpha_{q+1}-1}\left(1-t\right)^{\beta_{q}-\alpha_{q+1}-1}\Theta\left(\kappa_{l};-\frac{b}{t}-\frac{d}{1-t}\right)\\
 \times\sum_{m=0}^{\infty}\left(\alpha_{1}\right)_{m}\prod_{j=1}^{q-1}\frac{\mathcal{B}_{b,d}^{\left(\kappa_{l}\right)}\left(\alpha_{j+1}+m,\beta_{j}-\alpha_{j+1}\right)}{B\left(\alpha_{j+1},\beta_{j}-\alpha_{j+1}\right)}\frac{\left(zt\right)^{m}}{m!}dt\\
=\frac{\Gamma\left(\beta_{q}\right)}{\Gamma\left(\alpha_{q+1}\right)\Gamma\left(\beta_{q}-\alpha_{q+1}\right)}\int_{0}^{1}t^{\alpha_{q+1}-1}\left(1-t\right)^{\beta_{q}-\alpha_{q+1}-1}\Theta\left(\kappa_{l};-\frac{b}{t}-\frac{d}{1-t}\right)\\
\times{}_{q}F_{q-1}^{\left(\kappa_{l}\right)}\left[\begin{array}{ccc}
\alpha_{1} & \cdots & \alpha_{q}\\
\beta_{1} & \cdots & \beta_{q-1}
\end{array};zt;b,d\right]dt.
\end{multline}
For $p=q$, we have
\begin{align*}
_{q}F_{q}^{\left(\kappa_{l}\right)}\left[\begin{array}{ccc}
\alpha_{1} & \cdots & \alpha_{q}\\
\beta_{1} & \cdots & \beta_{q}
\end{array};z;b,d\right] & =\frac{\Gamma\left(\beta_{q}\right)}{\Gamma\left(\alpha_{q}\right)\Gamma\left(\beta_{q}-\alpha_{q}\right)}\int_{0}^{1}t^{\alpha_{q}-1}\left(1-t\right)^{\beta_{q}-\alpha_{q}-1}\Theta\left(\kappa_{l};-\frac{b}{t}-\frac{d}{1-t}\right)\\
 & \times{}_{q-1}F_{q-1}^{\left(\kappa_{l}\right)}\left[\begin{array}{ccc}
\alpha_{1} & \cdots & \alpha_{q-1}\\
\beta_{1} & \cdots & \beta_{q-1}
\end{array};zt;b,d\right]dt.
\end{align*}
It is clear that this relation is valid for $p<q$. This completes the proof.
\end{Proof}

\begin{Remark}
A multidimensional case of Euler type integral representation of (20) is given by:
\begin{multline}
_{q+1}F_{q}^{\left(\kappa_{l}\right)}\left[\begin{array}{ccc}
\alpha_{1} & \cdots & \alpha_{q+1}\\
\beta_{1} & \cdots & \beta_{q}
\end{array};z;b,d\right]
=\prod_{j=1}^{q}\frac{\Gamma \left(\beta_j\right)}{\Gamma \left(\alpha_{j+1}\right) \Gamma \left(\beta_j-\alpha_{j+1}\right)}\\
\times \int_0^1\cdots\int_0^1\prod_{j=1}^{q} \left\{t_{j}^{\alpha_{j+1}} (1-t_j)^{\beta_j-\alpha_{j+1}-1}\Theta\left(\kappa_l;-\frac{b}{t_j}-\frac{d}{1-t_j}\right) \right\}
\left(1-t_1 t_2 \cdots t_q z\right)^{-\alpha_1}dt_1 \cdots dt_q,\nonumber
\end{multline}
which follows from the repeated application of the functional equation (20). If we set $b=d=0$, then above equation reduces to the form given by \cite[p.132, Eq.(4.2)]{Extensions of certain classical integrals of Erdelyi}
\end{Remark}
\begin{Theorem}
For $p\leq q+1$, we have the following differentiation formula:
\begin{equation}
\frac{d^{n}}{dz^{n}}\left\{ _{p}F_{q}^{\left(\kappa_{l}\right)}\left[\begin{array}{ccc}
\alpha_{1} & \cdots & \alpha_{p}\\
\beta_{1} & \cdots & \beta_{q}
\end{array};z;b,d\right]\right\} =\frac{\left(\alpha_{1}\right)_{n}\cdots\left(\alpha_{p}\right)_{n}}{\left(\beta_{1}\right)_{n}\cdots\left(\beta_{q}\right)_{n}}{}_{p}F_{q}^{\left(\kappa_{l}\right)}\left[\begin{array}{ccc}
\alpha_{1}+n & \cdots & \alpha_{p}+n\\
\beta_{1}+n & \cdots & \beta_{q}+n
\end{array};z;b,d\right]
\end{equation}
\end{Theorem}
\begin{Proof} Taking the derivative of $_{q+1}F_{q}^{\left(\kappa_{l}\right)}$with respect to $z$, we obtain
\[
\begin{aligned}\frac{d}{dz}\left\{ _{q+1}F_{q}^{\left(\kappa_{l}\right)}\left[\begin{array}{ccc}
\alpha_{1} & \cdots & \alpha_{p}\\
\beta_{1} & \cdots & \beta_{q}
\end{array};z;b,d\right]\right\}  & =\frac{d}{dz}\left\{ \sum_{m=0}^{\infty}\left(\alpha_{1}\right)_{m}\prod_{j=1}^{q}\frac{\mathcal{B}_{b,d}^{\left(\kappa_{l}\right)}\left(\alpha_{j+1}+m,\beta_{j}-\alpha_{j+1}\right)}{B\left(\alpha_{j+1},\beta_{j}-\alpha_{j+1}\right)}\frac{z^{m}}{m!}\right\} \\
 & =\sum_{m=1}^{\infty}\left(\alpha_{1}\right)_{m}\prod_{j=1}^{q}\frac{\mathcal{B}_{b,d}^{\left(\kappa_{l}\right)}\left(\alpha_{j+1}+m,\beta_{j}-\alpha_{j+1}\right)}{B\left(\alpha_{j+1},\beta_{j}-\alpha_{j+1}\right)}\frac{z^{m-1}}{\left(m-1\right)!}.
\end{aligned}
\]
Replacing $m\rightarrow m+1$ we get
\begin{align*}
\frac{d}{dz}\left\{ _{q+1}F_{q}^{\left(\kappa_{l}\right)}\left[\begin{array}{ccc}
\alpha_{1} & \cdots & \alpha_{q+1}\\
\beta_{1} & \cdots & \beta_{q}
\end{array};z;b,d\right]\right\}  & =\sum_{m=0}^{\infty}\left(\alpha_{1}\right)_{m+1}\prod_{j=1}^{q}\frac{\mathcal{B}_{b,d}^{\left(\kappa_{l}\right)}\left(\alpha_{j+1}+m+1,\beta_{j}-\alpha_{j+1}\right)}{B\left(\alpha_{j+1},\beta_{j}-\alpha_{j+1}\right)}\frac{z^{m}}{m!}\\
 & =\alpha_{1}\frac{\prod_{j=1}^{q}\alpha_{j+1}}{\prod_{j=1}^{q}\beta_{j}}\sum_{m=0}^{\infty}\left(\alpha_{1}+1\right)\prod_{j=1}^{q}\frac{\mathcal{B}_{b,d}^{\left(\kappa_{l}\right)}\left(\alpha_{j+1}+1+m,\beta_{j}-\alpha_{j+1}\right)}{B\left(\alpha_{j+1}+1,\beta_{j}-\alpha_{j+1}\right)}\frac{z^{m}}{m!}\\
 & =\alpha_{1}\frac{\prod_{j=1}^{q}\alpha_{j+1}}{\prod_{j=1}^{q}\beta_{j}}{}_{q+1}F_{q}^{\left(\kappa_{l}\right)}\left[\begin{array}{ccc}
\alpha_{1}+1 & \cdots & \alpha_{q+1}+1\\
\beta_{1}+1 & \cdots & \beta_{q}+1
\end{array};z;b,d\right].
\end{align*}

Recursive application of this procedure gives us the general form
\[
\frac{d^{n}}{dz^{n}}\left\{ _{q+1}F_{q}^{\left(\kappa_{l}\right)}\left[\begin{array}{ccc}
\alpha_{1} & \cdots & \alpha_{q+1}\\
\beta_{1} & \cdots & \beta_{q}
\end{array};z;b,d\right]\right\} =\frac{\left(\alpha_{1}\right)_{n}\cdots\left(\alpha_{q+1}\right)_{n}}{\left(\beta_{1}\right)_{n}\cdots\left(\beta_{q}\right)_{n}}{}_{q+1}F_{q}^{\left(\kappa_{l}\right)}\left[\begin{array}{ccc}
\alpha_{1}+n & \cdots & \alpha_{q+1}+n\\
\beta_{1}+n & \cdots & \beta_{q}+n
\end{array};z;b,d\right].
\]

Similarly, we can prove this result for the case $p\leq q$.
\end{Proof}

For $p=2$ and $q=1$ we get the following corollary.
\newtheorem{corollary}{Corollary}
\begin{corollary}
\begin{equation}
\frac{d^{n}}{dz^{n}}\left\{ _{2}F_{1}^{\left(\kappa_{l}\right)}\left[\begin{array}{c}
\alpha_{1},\alpha_{2}\\
\beta_{1}
\end{array};z;b,d\right]\right\} =\frac{\left(\alpha_{1}\right)_{n}\left(\alpha_{2}\right)_{n}}{\left(\beta_{1}\right)_{n}}{}_{2}F_{1}^{\left(\kappa_{l}\right)}\left[\begin{array}{c}
\alpha_{1}+n,\alpha_{2}+n\\
\beta_{1}+n
\end{array};z;b,d\right]
\end{equation}
\end{corollary}
\begin{Remark}
If we take $\Theta\left(\kappa_{l};z\right)={}_{1}F_{1}\left(\alpha;\beta;z\right),b=d$,
then equation (21) reduces to
\[
\frac{d^{n}}{dz^{n}}\left\{ F_{b}^{\left(\alpha,\beta\right)}\left[\begin{array}{c}
a,b\\
c
\end{array};z\right]\right\} =\frac{\left(b\right)_{n}\left(a\right)_{n}}{\left(c\right)_{n}}F_{b}^{\left(\alpha,\beta\right)}\left[\begin{array}{c}
a+n,b+n\\
c+n
\end{array};z\right]
\]
given in \cite[Theorem 3.3]{Extension of gamma beta and hypergeometric functions}.
\end{Remark}

In order to derive the Mellin-Barnes type contour integral representation of $_{p}F_{q}^{\left(\kappa_{l}\right)}\left(\alpha_{1},\cdots,\alpha_{p};\beta_{1},\cdots,\beta_{q};z;b,d\right)$
we need to introduce the following well-known theorem. It is widely used
to evaluate definite integrals and infinite series.
\begin{Theorem}[\textbf{Ramanujan's Master Theorem} \cite{Ramanujan Master Theorem}]\emph{Assume $f$ admits an expansion
of the form: }
\[
f\left(x\right)=\sum_{k=0}^{\infty}\frac{\lambda\left(k\right)}{k!}\left(-x\right)^{k}.
\]
\emph{Then, the Mellin transform of $f$ is given by}
\[
\int_{0}^{\infty}x^{s-1}f\left(x\right)dx=\Gamma\left(s\right)\lambda\left(-s\right).
\]
\end{Theorem}
By means of Ramanujan's master theorem we obtain:
\begin{Theorem}
We have the following Mellin-Barnes type integral representation
of function (18), namely,
\begin{multline}
_{p}F_{q}^{\left(\kappa_{l}\right)}\left[\begin{array}{ccc}
\alpha_{1} & \cdots & \alpha_{p}\\
\beta_{1} & \cdots & \beta_{q}
\end{array};z;b,d\right]\\
=\begin{cases}
\displaystyle\frac{1}{2\pi i}\int_{L}\prod_{j=1}^{q}\frac{\mathcal{B}_{b,d}^{\left(\kappa_{l}\right)}\left(\alpha_{j+1}-s,\beta_{j}-\alpha_{j+1}\right)}{B\left(\alpha_{j+1},\beta_{j}-\alpha_{j+1}\right)}\frac{\Gamma\left(s\right)\Gamma\left(\alpha_{1}-s\right)}{\Gamma\left(\alpha_{1}\right)}\left(-z\right)^{-s}ds, & p=q+1\\[15pt]
\displaystyle\frac{1}{2\pi i}\int_{L}\prod_{j=1}^{q}\frac{\mathcal{B}_{b,d}^{\left(\kappa_{l}\right)}\left(\alpha_{j}-s,\beta_{j}-\alpha_{j}\right)}{B\left(\alpha_{j},\beta_{j}-\alpha_{j}\right)}\Gamma\left(s\right)\left(-z\right)^{-s}ds, & p=q\\[15pt]
\displaystyle\frac{1}{2\pi i}\int_{L}\prod_{j=1}^{q}\frac{\mathcal{B}_{b,d}^{\left(\kappa_{l}\right)}\left(\alpha_{j}-s,\beta_{j+r}-\alpha_{j}\right)}{B\left(\alpha_{j},\beta_{j+r}-\alpha_{j}\right)}\frac{\Gamma\left(\beta_{1}\right)\cdots\Gamma\left(\beta_{r}\right)\Gamma\left(s\right)}{\Gamma\left(\beta_{1}-s\right)\cdots\Gamma\left(\beta_{r}-s\right)}\left(-z\right)^{-s}ds, & r=q-p,\ p<q.
\end{cases}
\end{multline}
where $L$ is a Barnes path of integration, that is, $L$ starts at $-i\infty$ and runs to $+i\infty$ in the $s$-plane, curving if necessary to put the poles of $\Gamma(\alpha_1-s)$ to the left of the path and to put the poles of $\Gamma(s)$ to the right of the path.
\end{Theorem}
\begin{Proof}The result follows directly by using Ramanujan's Master Theorem
and inverse Mellin transform.
\end{Proof}

In section 3, a Multidimensional Ramanjuan's Master Theorem will be
introduced to obtain double and multiple contour integral representations
for extended Appell's and Lauricella's hypergeometric functions.

\subsection{{\normalsize \textbf{More properties about Gauss hypergeometric function}}}
Although many different kinds of extensions of Gauss hypergeometric function have been given, we still know little about their sorts of properties. Therefore, we list some formulas and conclusions about extended Gauss hypergeometric function,of which some are originated from classical theories and others are drawn on the basis of recent researches.
\begin{Theorem}
The following transformations hold true for extended Gauss hypergeometric function:
\begin{description}
\item [1. Pfaff transformation]
\begin{equation}
_{2}F_{1}^{\left(\kappa_{l}\right)}\left[\begin{array}{c}
\alpha_{1},\alpha_{2}\\
\beta_{1}
\end{array};z;b,d\right]=\left(1-z\right)^{-\alpha_{1}}{}_{2}F_{1}^{\left(\kappa_{l}\right)}\left[\begin{array}{c}
\alpha_{1},\beta_{1}-\alpha_{2}\\
\alpha_{2}
\end{array};\frac{z}{1-z};d,b\right],\ \left(\left|\arg\left(1-z\right)\right|<\pi\right).
\end{equation}
\item [2. Euler transformation] When $\Theta\left(\kappa_{l};z\right)=\exp\left(z\right)$,
we have:
\begin{equation}
_{2}F_{1}\left[\begin{array}{c}
\alpha_{1},\alpha_{2}\\
\beta_{1}
\end{array};z;b,d\right]=e^{-\left(1-z\right)b-zd}\left(1-z\right)^{\beta_{1}-\alpha_{2}-\alpha_{1}}{}_{2}F_{1}\left[\begin{array}{c}
\beta_{1}-\alpha_{1},\beta_{1}-\alpha_{2}\\
\beta_{1}
\end{array};z;\frac{b}{1-z},\left(1-z\right)d\right].
\end{equation}
\[
\left(\left|\arg\left(1-z\right)\right|<\pi,\ \Re\left(z\right)<1\right)
\]
\end{description}
\end{Theorem}
\begin{Proof}
1. The proof of transformation (24) is direct. Since we have $\left[1-z\left(1-t\right)\right]^{-a}=\left(1-z\right)^{-a}\left(1+\frac{z}{1-z}t\right)^{-a}$
then, by replacing $t\rightarrow1-t$, we obtain
\begin{align*}
_{2}F_{1}^{\left(\kappa_{l}\right)}\left[\begin{array}{c}
\alpha_{1},\alpha_{2}\\
\beta_{1}
\end{array};z;b,d\right] & =\frac{1}{B\left(\alpha_{2},\beta_{1}-\alpha_{2}\right)}\int_{0}^{1}\left(1-t\right)^{\alpha_{2}-1}t^{\beta_{1}-\alpha_{2}-1}\left[1-z\left(1-t\right)\right]^{-\alpha_{1}}\Theta\left(\kappa_{l};-\frac{b}{1-t}-\frac{d}{t}\right)dt\\
 & =\frac{\left(1-z\right)^{-a}}{B\left(\alpha_{2},\beta_{1}-\alpha_{2}\right)}\int_{0}^{1}t^{\beta_{1}-\alpha_{2}-1}\left(1-t\right)^{\alpha_{2}-1}\left(1+\frac{z}{1-z}t\right)^{-\alpha_{1}}\Theta\left(\kappa_{l};-\frac{d}{t}-\frac{b}{1-t}\right)dt\\
 & =\left(1-z\right)^{-a}{}_{2}F_{1}^{\left(\kappa_{l}\right)}\left[\begin{array}{c}
\alpha_{1},\beta_{1}-\alpha_{2}\\
\alpha_{2}
\end{array};\frac{z}{1-z};d,b\right].
\end{align*}

2. When $\Theta\left(\kappa_{l};z\right)=\exp\left(z\right)$, from
Pfaff transformation (24) we have
\[
{}_{2}F_{1}\left[\begin{array}{c}
\alpha_{1},\alpha_{2}\\
\beta_{1}
\end{array};z;b,d\right]=\left(1-z\right)^{-\alpha_{1}}{}_{2}F_{1}\left[\begin{array}{c}
\alpha_{1},\beta_{1}-\alpha_{2}\\
\alpha_{2}
\end{array};\frac{z}{1-z};d,b\right].
\]
The rest of the proof will be presented in the \textbf{section 4} since we need
more extra results built in that section.
\end{Proof}

We also give a somewhat complicated differentiation formula for extended
Gauss hypergeometric function.
\begin{Theorem}
One has
\begin{equation}
\frac{d^{n}}{dz^{n}}\left\{ z^{\alpha_{1}+n-1}{}_{2}F_{1}^{\left(\kappa_{l}\right)}\left[\begin{array}{c}
\alpha_{1},\alpha_{2}\\
\beta_{1}
\end{array};z;b,d\right]\right\} =\left(\alpha_{1}\right)_{n}z^{\alpha_{1}-1}{}_{2}F_{1}^{\left(\kappa_{l}\right)}\left[\begin{array}{c}
\alpha_{1},\alpha_{2}\\
\beta_{1}
\end{array};z;b,d\right].
\end{equation}
\end{Theorem}

\begin{Proof}This formula can be obtained by direct computation,
\begin{align*}
\frac{d}{dz}\left\{ z^{\alpha_{1}}{}_{2}F_{1}^{\left(\kappa_{l}\right)}\left[\begin{array}{c}
\alpha_{1},\alpha_{2}\\
\beta_{1}
\end{array};z;b,d\right]\right\}  & =\frac{d}{dz}\left\{ \sum_{n=0}^{\infty}\left(\alpha_{1}\right)_{n}\frac{\mathcal{B}_{b,d}^{\left(\kappa_{l}\right)}\left(\alpha_{2}+n,\beta_{1}-\alpha_{2}\right)}{B\left(\alpha_{2},\beta_{1}-\alpha_{2}\right)}\frac{z^{n+\alpha_{1}}}{n!}\right\} \\
 & =\sum_{n=0}^{\infty}\left(\alpha_{1}\right)_{n}\left(\alpha_{1}+n\right)\frac{\mathcal{B}_{b,d}^{\left(\kappa_{l}\right)}\left(\alpha_{2}+n,\beta_{1}-\alpha_{2}\right)}{B\left(\alpha_{2},\beta_{1}-\alpha_{2}\right)}\frac{z^{n+\alpha_{1}-1}}{n!}\\
\left(Note\ that\ \left(a+n\right)\left(a\right)_{n}=a\left(a+1\right)_{n}\right) & =\alpha_{1}z^{\alpha_{1}-1}\sum_{n=0}^{\infty}\left(\alpha_{1}+1\right)_{n}\frac{\mathcal{B}_{b,d}^{\left(\kappa_{l}\right)}\left(\alpha_{2}+n,\beta_{1}-\alpha_{2}\right)}{B\left(\alpha_{2},\beta_{1}-\alpha_{2}\right)}\frac{z^{n}}{n!}\\
 & =\alpha_{1}z^{\alpha_{1}-1}{}_{2}F_{1}^{\left(\kappa_{l}\right)}\left[\begin{array}{c}
\alpha_{1}+1,\alpha_{2}\\
\beta_{1}
\end{array};z;b,d\right].
\end{align*}

\begin{align*}
\frac{d^{2}}{dz^{2}}\left\{ z^{\alpha_{1}+1}{}_{2}F_{1}^{\left(\kappa_{l}\right)}\left[\begin{array}{c}
\alpha_{1},\alpha_{2}\\
\beta_{1}
\end{array};z;b,d\right]\right\}  & =\sum_{n=0}^{\infty}\left(\alpha_{1}\right)_{n}\left(\alpha_{1}+n+1\right)\left(\alpha_{1}+n\right)\frac{\mathcal{B}_{b,d}^{\left(\kappa_{l}\right)}\left(\alpha_{2}+n,\beta_{1}-\alpha_{2}\right)}{B\left(\alpha_{2},\beta_{1}-\alpha_{2}\right)}\frac{z^{n+\alpha_{1}-1}}{n!}\\
 & =\alpha_{1}\left(\alpha_{1}+1\right)z^{\alpha_{1}}{}_{2}F_{1}^{\left(\kappa_{l}\right)}\left[\begin{array}{c}
\alpha_{1}+2,\alpha_{2}\\
\beta_{1}
\end{array};z;b,d\right]
\end{align*}
In the procedure, we have used the fact that $\left(a+n+1\right)\left(a+n\right)\left(a\right)_{n}=a\left(a+1\right)\left(a+2\right)_{n}$.
Recursive application of this procedure give us the general form (26).
\end{Proof}

A class of recurrence relations for Gauss hypergeometric functions
\[
_{2}F_{1}\left[\begin{array}{c}
\alpha_{1}\pm n,\alpha_{2}\\
\beta_{1}
\end{array};z\right],\ {}_{2}F_{1}\left[\begin{array}{c}
\alpha_{1},\alpha_{2}\pm n\\
\beta_{1}
\end{array};z\right]\ and\ {}_{2}F_{1}\left[\begin{array}{c}
\alpha_{1},\alpha_{2}\\
\beta_{1}\pm n
\end{array};z\right]\ \ \left(n\in\mathbb{N}_{0}\right)
\]
have been established in \cite{Recursion formulas for Appells hypergeometric function with some applications to radiation field problems} as a useful tool to find some new recursion formulas for Appell hypergeometric functions
\[
\mathbf{F}_{2}\left(\sigma,\alpha_{1},\alpha_{2}\pm n;\beta_{1},\beta_{2};x,y\right)\ and\ \mathbf{F}_{2}\left(\sigma,\alpha_{1},\alpha_{2};\beta_{1},\beta_{2}\pm n;x,y\right).
\]

It is natural to consider whether these results are true for our new
definitions. The following theorem gives us a positive answer.
\begin{Theorem}
The following recurrence relations hold true for the
extended Gauss hypergeometric functions $_{2}F_{1}^{\left( \kappa_{l}\right)}$:
\begin{description}
\item [1. Recurrence relations for $_{2}F_{1}^{\left(\kappa_{l}\right)}\left(\alpha_{1}\pm n,\alpha_{2};\beta_{1};z;b,d\right)$]
\begin{equation}
_{2}F_{1}^{\left(\kappa_{l}\right)}\left[\begin{array}{c}
\alpha_{1}+n,\alpha_{2}\\
\beta_{1}
\end{array};z;b,d\right]={}_{2}F_{1}^{\left(\kappa_{l}\right)}\left[\begin{array}{c}
\alpha_{1},\alpha_{2}\\
\beta_{1}
\end{array};z;b,d\right]+\frac{\alpha_{2}z}{\beta_{1}}\sum_{k=1}^{n}{}_{2}F_{1}^{\left(\kappa_{l}\right)}\left[\begin{array}{c}
\alpha_{1}+n-k+1,\alpha_{2}+1\\
\beta_{1}+1
\end{array};z;b,d\right],
\end{equation}
\begin{equation}
_{2}F_{1}^{\left(\kappa_{l}\right)}\left[\begin{array}{c}
\alpha_{1}-n,\alpha_{2}\\
\beta_{1}
\end{array};z;b,d\right]={}_{2}F_{1}^{\left(\kappa_{l}\right)}\left[\begin{array}{c}
\alpha_{1},\alpha_{2}\\
\beta_{1}
\end{array};z;b,d\right]-\frac{\alpha_{2}z}{\beta_{1}}\sum_{k=1}^{n}{}_{2}F_{1}^{\left(\kappa_{l}\right)}\left[\begin{array}{c}
\alpha_{1}-k+1,\alpha_{2}+1\\
\beta_{1}+1
\end{array};z;b,d\right].
\end{equation}
\[
\left(\left|z\right|<1;\ n\in\mathbb{N}_{0}\right)
\]
\item [2. Recurrence relation for $_{2}F_{1}^{\left(\kappa_{l}\right)}\left(\alpha_{1},\alpha_{2};\beta_{1}+n;z;b,d\right)$]
\begin{equation}
_{2}F_{1}^{\left(\kappa_{l}\right)}\left[\begin{array}{c}
\alpha_{1},\alpha_{2}\\
\beta_{1}+n
\end{array};z;b,d\right]=\frac{\left(\beta_{1}\right)_{n}}{\left(\beta_{1}-\alpha_{2}\right)_{n}}\sum_{k=0}^{n}\left(-1\right)^{k}\left(\begin{array}{c}
n\\
k
\end{array}\right)\frac{\left(\alpha_{2}\right)_{k}}{\left(\beta_{1}\right)_{k}}{}_{2}F_{1}^{\left(\kappa_{l}\right)}\left[\begin{array}{c}
\alpha_{1},\alpha_{2}+k\\
\beta_{1}+k
\end{array};z;b,d\right].
\end{equation}
\[
\left(\left|z\right|<1;\ n\in\mathbb{N}_{0}\right)
\]
\item [3. Recurrence relation for $_{2}F_{1}^{\left(\kappa_{l}\right)}\left(\alpha_{1},\alpha_{2}+n;\beta_{1};z;b,d\right)$]
\begin{equation}
_{2}F_{1}^{\left(\kappa_{l}\right)}\left[\begin{array}{c}
\alpha_{1},\alpha_{2}+n\\
\beta_{1}
\end{array};z;b,d\right]=\frac{\left(\beta_{1}-\alpha_{2}\right)_{2n}}{\left(\beta_{1}-\alpha_{2}\right)_{n}\left(\alpha_{2}\right)_{n}}\sum_{i=1}^{n}\left(-n\right)_{i}\frac{\left(\alpha_{2}\right)_{i+n}}{\left(\beta_{1}\right)_{i+n}i!}{}_{2}F_{1}^{\left(\kappa_{l}\right)}\left[\begin{array}{c}
\alpha_{1},\alpha_{2}+i+n\\
\beta_{1}+i+n
\end{array};z;b,d\right].
\end{equation}
\[
\left(\left|z\right|<1;\ n\in\mathbb{N}_{0}\right)
\]
\end{description}
\end{Theorem}
\begin{Proof}
1. By means of Euler integral representation of the extended Gauss
hypergeometric functions given by (13), we have, by writing $\left(1-zt\right)^{-\alpha_{1}-n}=\left(1-zt\right)^{-\alpha_{1}-n-1}-zt\left(1-zt\right)^{-\alpha_{1}-n-1}$,
\begin{align*}
_{2}F_{1}^{\left(\kappa_{l}\right)}\left[\begin{array}{c}
\alpha_{1}+n,\alpha_{2}\\
\beta_{1}
\end{array};z;b,d\right] & =\frac{\Gamma\left(\beta_{1}\right)}{\Gamma\left(\alpha_{2}\right)\Gamma\left(\beta_{1}-\alpha_{2}\right)}\int_{0}^{1}t^{\alpha_{2}-1}\left(1-t\right)^{\beta_{1}-\alpha_{2}-1}\left(1-zt\right)^{-\alpha_{1}-n-1}\Theta\left(\kappa_{l};-\frac{b}{t}-\frac{d}{1-t}\right)dt\\
 & \ -\frac{z\Gamma\left(\beta_{1}\right)}{\Gamma\left(\alpha_{2}\right)\Gamma\left(\beta_{1}-\alpha_{2}\right)}\int_{0}^{1}t^{\alpha_{2}}\left(1-t\right)^{\beta_{1}-\alpha_{2}-1}\left(1-zt\right)^{-\alpha_{1}-n-1}\Theta\left(\kappa_{l};-\frac{b}{t}-\frac{d}{1-t}\right)dt\\
 & ={}_{2}F_{1}^{\left(\kappa_{l}\right)}\left(\begin{array}{c}
\alpha_{1}+n+1,\alpha_{2}\\
\beta_{1}
\end{array};z;b,d\right)-\left(\frac{z}{\alpha_{1}+n}\right)\frac{d}{dz}\left\{ _{2}F_{1}^{\left(\kappa_{l}\right)}\left[\begin{array}{c}
\alpha_{1}+n,\alpha_{2}\\
\beta_{1}
\end{array};z;b,d\right]\right\} .
\end{align*}
Replacing $n$ by $n-1$, we find that
\begin{multline}
_{2}F_{1}^{\left(\kappa_{l}\right)}\left[\begin{array}{c}
\alpha_{1}+n,\alpha_{2}\\
\beta_{1}
\end{array};z;b,d\right]=
{}_{2}F_{1}^{\left(\kappa_{l}\right)}\left[\begin{array}{c}
\alpha_{1}+n-1,\alpha_{2}\\
\beta_{1}
\end{array};z;b,d\right]\\+\left(\frac{z}{\alpha_{1}+n-1}\right)\frac{d}{dz}\left\{ _{2}F_{1}^{\left(\kappa_{l}\right)}\left[\begin{array}{c}
\alpha_{1}+n-1,\alpha_{2}\\
\beta_{1}
\end{array};z;b,d\right]\right\} .
\end{multline}
Applying this identity recursively, we obtain:
\begin{multline}
_{2}F_{1}^{\left(\kappa_{l}\right)}\left[\begin{array}{c}
\alpha_{1}+n,\alpha_{2}\\
\beta_{1}
\end{array};z;b,d\right]={}_{2}F_{1}^{\left(\kappa_{l}\right)}\left[\begin{array}{c}
\alpha_{1},\alpha_{2}\\
\beta_{1}
\end{array};z;b,d\right]\\+z\sum_{k=1}^{n}\frac{1}{\alpha_{1}+n-k}\frac{d}{dz}\left\{ _{2}F_{1}^{\left(\kappa_{l}\right)}\left[\begin{array}{c}
\alpha_{1}+n-k,\alpha_{2}\\
\beta_{1}
\end{array};z;b,d\right]\right\} \\
={}_{2}F_{1}^{\left(\kappa_{l}\right)}\left[\begin{array}{c}
\alpha_{1},\alpha_{2}\\
\beta_{1}
\end{array};z;b,d\right]+z\sum_{k=1}^{n}\frac{\alpha_{2}}{\beta_{1}}{}_{2}F_{1}^{\left(\kappa_{l}\right)}\left[\begin{array}{c}
\alpha_{1}+n-k+1,\alpha_{2}+1\\
\beta_{1}+1
\end{array};z;b,d\right].
\end{multline}
To prove the second formula, we still begin with the integral representation of extended Gauss hypergeometric function. Writing $\left(1-zt\right)^{-\left(\alpha_{1}-n\right)}=\left(1-zt\right)^{-\left(\alpha_{1}-n+1\right)}-zt\left(1-zt\right)^{-\left(\alpha_{1}-n+1\right)}$
we have
\begin{multline}
_{2}F_{1}^{\left(\kappa_{l}\right)}\left[\begin{array}{c}
\alpha_{1}-n,\alpha_{2}\\
\beta_{1}
\end{array};z;b,d\right]=\\
\frac{\Gamma\left(\beta_{1}\right)}{\Gamma\left(\alpha_{2}\right)\Gamma\left(\beta_{1}-\alpha_{2}\right)}\int_{0}^{1}t^{\alpha_{2}-1}\left(1-t\right)^{\beta_{1}-\alpha_{2}-1}\left(1-zt\right)^{-\left(\alpha_{1}-n+1\right)}\Theta\left(\kappa_{l};-\frac{b}{t}-\frac{d}{1-t}\right)dt\\ -\frac{z\Gamma\left(\beta_{1}\right)}{\Gamma\left(\alpha_{2}\right)\Gamma\left(\beta_{1}-\alpha_{2}\right)}\int_{0}^{1}t^{\alpha_{2}}\left(1-t\right)^{\beta_{1}-\alpha_{2}-1}\left(1-zt\right)^{-\left(\alpha_{1}-n+1\right)}\Theta\left(\kappa_{l};-\frac{b}{t}-\frac{d}{1-t}\right)dt\\
={}_{2}F_{1}^{\left(\kappa_{l}\right)}\left[\begin{array}{c}
\alpha_{1}-n+1,\alpha_{2}\\
\beta_{1}
\end{array};z;b,d\right]-\left(\frac{z}{\alpha_{1}-n}\right)\frac{d}{dz}\left\{ _{2}F_{1}^{\left(\kappa_{l}\right)}\left[\begin{array}{c}
\alpha_{1}-n,\alpha_{2}\\
\beta_{1}
\end{array};z;b,d\right]\right\} .
\end{multline}
Replacing $n$ by $n+1$, we find that
\begin{multline}
_{2}F_{1}^{\left(\kappa_{l}\right)}\left[\begin{array}{c}
\alpha_{1}-n,\alpha_{2}\\
\beta_{1}
\end{array};z;b,d\right]={}_{2}F_{1}^{\left(\kappa_{l}\right)}\left[\begin{array}{c}
\alpha_{1}-n-1,\alpha_{2}\\
\beta_{1}
\end{array};z;b,d\right]\\
-z\frac{d}{dz}\left\{ \frac{1}{n-\alpha_{1}+1}{}_{2}F_{1}^{\left(\kappa_{l}\right)}\left[\begin{array}{c}
\alpha_{1}-n-1,\alpha_{2}\\
\beta_{1}
\end{array};z;b,d\right]\right\} .
\end{multline}
Repeating this recurrence relation $n$ times and applying \textbf{Corollary 1}, we obtain:
\begin{multline}
_{2}F_{1}^{\left(\kappa_{l}\right)}\left[\begin{array}{c}
\alpha_{1}-n,\alpha_{2}\\
\beta_{1}
\end{array};z;b,d\right]={}_{2}F_{1}^{\left(\kappa_{l}\right)}\left[\begin{array}{c}
\alpha_{1}-2n,\alpha_{2}\\
\beta_{1}
\end{array};z;b,d\right]\\
+\frac{\alpha_{2}z}{\beta_{1}}\sum_{k=1}^{n}{}_{2}F_{1}^{\left(\kappa_{l}\right)}\left[\begin{array}{c}
\alpha_{1}-n-k+1,\alpha_{2}+1\\
\beta_{1}+1
\end{array};z;b,d\right].
\end{multline}
The identity (28) follows directly from (35) upon replacing $\alpha_{1}$ by $\alpha_{1}+n$ $\left(n\in\mathbb{N}_{0}\right)$.\\

2. Using Euler integral representation (13) and the fact that
\[
\left(1-t\right)^{\beta_{1}+n-\alpha_{2}-1}=\left(1-t\right)^{\beta_{1}-\alpha_{1}-1}\left(1-t\right)^{n}=\left(1-t\right)^{\beta_{1}-\alpha_{2}-1}\sum_{k=0}^{n}\dbinom{n}{k}\left(-1\right)^{k}t^{k},
\]
we find
\begin{align*}
_{2}F_{1}^{\left(\kappa_{l}\right)}\left[\begin{array}{c}
\alpha_{1},\alpha_{2}\\
\beta_{1}+n
\end{array};z;b,d\right] & =\frac{\Gamma\left(\beta_{1}+n\right)}{\Gamma\left(\alpha_{2}\right)\Gamma\left(\beta_{1}-\alpha_{2}+n\right)}\sum_{k=0}^{n}\left(\begin{array}{c}
n\\
k
\end{array}\right)\left(-1\right)^{k}\int_{0}^{1}t^{\alpha_{1}+k-1}\left(1-t\right)^{\beta_{1}+k-\left(\alpha_{2}+k\right)-1}\\
 & \ \ \times\left(1-zt\right)^{\alpha-1}\Theta\left(\kappa_{l};-\frac{b}{t}-\frac{d}{1-t}\right)dt\\
 & =\frac{\Gamma\left(\beta_{1}+n\right)}{\Gamma\left(\alpha_{2}\right)\Gamma\left(\beta_{1}-\alpha_{2}+n\right)}\sum_{k=0}^{n}\left(\begin{array}{c}
n\\
k
\end{array}\right)\left(-1\right)^{k}\frac{\Gamma\left(\alpha_{2}+k\right)\Gamma\left(\beta_{1}-\alpha_{2}\right)}{\Gamma\left(\beta_{1}+k\right)}{}_{2}F_{1}^{\left(\kappa_{l}\right)}\left[\begin{array}{c}
\alpha_{1},\alpha_{2}+k\\
\beta_{1}+k
\end{array};z;b,d\right]\\
 & =\frac{\left(\beta_{1}\right)_{n}}{\left(\beta_{1}-\alpha_{2}\right)_{n}}\sum_{k=0}^{n}\left(\begin{array}{c}
n\\
k
\end{array}\right)\left(-1\right)^{k}\frac{\left(\alpha_{2}\right)_{k}}{\left(\beta_{1}\right)_{k}}{}_{2}F_{1}^{\left(\kappa_{l}\right)}\left[\begin{array}{c}
\alpha_{1},\alpha_{2}+k\\
\beta_{1}+k
\end{array};z;b,d\right].
\end{align*}
This completes the proof of Eq (29).

3. The demonstration of assertion 3 is based on the following elementary expansion
\begin{equation}
\left(1-t\right)^{-n}=\sum_{i=0}^{\infty}\dbinom{-n}{i}\left(-1\right)^{i}t^{i}=\sum_{i=0}^{n}\left(-n\right)_{i}\frac{t^{i}}{i!}.
\end{equation}
Substituting expansion (36) in the integral representation of extended Gauss hypergeometric function (13) we have
\begin{align*}
_{2}F_{1}^{\left(\kappa_{l}\right)}\left[\begin{array}{c}
\alpha_{1},\alpha_{2}+n\\
\beta_{1}
\end{array};z;b,d\right] & =\frac{\Gamma\left(\beta_{1}\right)}{\Gamma\left(\alpha_{2}+n\right)\Gamma\left(\beta_{1}-\alpha_{2}-n\right)}\int_{0}^{1}t^{\alpha_{2}+n-1}\left(1-t\right)^{\beta_{1}-\alpha_{2}-n-1}\left(1-zt\right)^{-\alpha_{1}}\\
 & \ \ \times\Theta\left(\kappa_{l};-\frac{b}{t}-\frac{d}{1-t}\right)dt\\
 &
 =\frac{\Gamma\left(\beta_{1}\right)}{\Gamma\left(\alpha_{2}\right)\left(\alpha_{2}\right)_{n}\Gamma\left(\beta_{1}-\alpha_{2}-n\right)}\sum_{i=0}^{n}\left(-n\right)_{i}\int_{0}^{1}t^{\alpha_{2}+n+i-1}\left(1-t\right)^{\beta_{1}-\alpha_{2}-1}\\
 & \ \ \times\left(1-zt\right)^{-\alpha_{1}}\Theta\left(\kappa_{l};-\frac{b}{t}-\frac{d}{1-t}\right)dt\\
 & =\frac{\Gamma\left(\beta_{1}\right)}{\Gamma\left(\alpha_{2}\right)\left(\alpha_{2}\right)_{n}\Gamma\left(\beta_{1}-\alpha_{2}-n\right)}\sum_{i=0}^{n}\left(-n\right)_{i}\frac{\Gamma\left(\alpha_{2}+i+n\right)\Gamma\left(\beta_{1}-\alpha_{2}\right)}{\Gamma\left(\beta_{1}+i+n\right)}\\
 & \ \ \times{}_{2}F_{1}^{\left(\kappa_{l}\right)}\left[\begin{array}{c}
\alpha_{1},\alpha_{2}+n+i\\
\beta_{1}+n+i
\end{array};z;b,d\right]\\
 & =\frac{\Gamma\left(\beta_{1}-\alpha_{2}\right)}{\left(\alpha_{2}\right)_{n}\Gamma\left(\beta_{1}-\alpha_{2}-n\right)}\sum_{i=0}^{n}\left(-n\right)_{i}\frac{\left(\alpha_{2}\right)_{i+n}}{\left(\beta_{1}\right)_{i+n}}{}_{2}F_{1}^{\left(\kappa_{l}\right)}\left[\begin{array}{c}
\alpha_{1},\alpha_{2}+n+i\\
\beta_{1}+n+i
\end{array};z;b,d\right].
\end{align*}
And equation (30) follows by using the following property
\[
\frac{\Gamma\left(\beta_{1}-\alpha_{2}\right)}{\Gamma\left(\beta_{1}-\alpha_{2}-n\right)}=\left(\beta_{1}-\alpha_{2}+n\right)_{n}=\frac{\left(\beta_{1}-\alpha_{2}\right)_{2n}}{\left(\beta_{1}-\alpha_{2}\right)_{n}}.
\]
\end{Proof}
\begin{Remark}
For classical Gauss hypergeometric function, the parameters
$\alpha_{1}$ and $\alpha_{2}$ are symmetric, i.e.,
\[
{}_{2}F_{1}\left(\alpha_{1},\alpha_{2};\beta_{1};z\right)={}_{2}F_{1}\left(\alpha_{2},\alpha_{1};\beta_{1};z\right).
\]
Thus we only need to establish the recurrence relations for one of them and the results can be valid for both $\alpha_{1}$ and $\alpha_{2}$. However, the numerator parameters for extended Gauss hypergeometric functions possess no such symmetry.
\end{Remark}

\subsection{{\normalsize \textbf{A further generalization and its connection with fractional
calculus}}}
Our definition of extended generalized hypergeometric function can in fact be generalized to the following form:
\begin{Definition}
For suitably constrained (real or complex) parameters
$\alpha_{j},\ j=1,\cdots,p$; $\beta_{i},\ i=1,\cdots,q$, we define
the extended generalized hypergeometric functions by
\begin{multline}
_{p}F_{q}^{\left(\kappa_{l}\right)}\left[\begin{array}{ccc}
\left(\alpha_{1},k_{1}\right) & \cdots & \left(\alpha_{p},k_{p}\right)\\
\beta_{1}, & \cdots & \beta_{q},
\end{array};z;b,d\right]=\\
\begin{cases}
\displaystyle\sum_{m=0}^{\infty}\left(\alpha_{1}\right)_{k_{1}m}\prod_{j=1}^{q}\frac{\mathcal{B}_{b,d}^{\left(\kappa_{l}\right)}\left(\alpha_{j+1}+k_{j+1}m,\beta_{j}-\alpha_{j+1}\right)}{B\left(\alpha_{j+1},\beta_{j}-\alpha_{j+1}\right)}\frac{z^{m}}{m!}, & \left|z\right|<1;\ p=q+1\\[15pt]
\displaystyle\sum_{m=0}^{\infty}\prod_{j=1}^{q}\frac{\mathcal{B}_{b,d}^{\left(\kappa_{l}\right)}\left(\alpha_{j}+k_{j}m,\beta_{j}-\alpha_{j}\right)}{B\left(\alpha_{j},\beta_{j}-\alpha_{j}\right)}\frac{z^{m}}{m!}, & p=q\\[15pt]
\displaystyle\sum_{m=0}^{\infty}\frac{1}{\left(\beta_{1}\right)_{m}}\cdots\frac{1}{\left(\beta_{r}\right)_{m}}\prod_{j=1}^{p}\frac{\mathcal{B}_{b,d}^{\left(\kappa_{l}\right)}\left(\alpha_{j}+k_{j}m,\beta_{r+j}-\alpha_{j}\right)}{B\left(\alpha_{j},\beta_{r+j}-\alpha_{j}\right)}\frac{z^{m}}{m!}, & r=q-p,\ p<q
\end{cases}
\end{multline}
where the new introduced parameters $k_{j},\ j=1,\cdots,p$ are non-negative integers.
\end{Definition}

Obviously, (37) reduces to (18) whenever $k_{j}=1,\ j=1,\cdots,p$.
To illustrate its advantages we first consider the following function:
\begin{equation}
_{2}F_{1}^{\left(\kappa_{l}\right)}\left[\begin{array}{c}
\left(\alpha_{1},1\right),\left(\alpha_{2},k_{2}\right)\\
\beta_{1}
\end{array};z;b,d\right]=\sum_{n=0}^{\infty}\left(\alpha_{1}\right)_{n}\frac{\mathcal{B}_{b,d}^{\left(\kappa_{l}\right)}\left(\alpha_{2}+k_{2}n,\beta_{1}-\alpha_{2}\right)}{B\left(\alpha_{2},\beta_{1}-\alpha_{2}\right)}\frac{z^{n}}{n!}.
\end{equation}
Its integral representation can be written as
\begin{multline}
_{2}F_{1}^{\left(\kappa_{l}\right)}\left[\begin{array}{c}
\left(\alpha_{1},1\right),\left(\alpha_{2},k_{2}\right)\\
\beta_{1}
\end{array};z;b,d\right]=\frac{1}{B\left(\alpha_{2},\beta_{1}-\alpha_{2}\right)}\\
\int_{0}^{1}t^{\alpha_{2}-1}\left(1-t\right)^{\beta_{1}-\alpha_{2}-1}\left(1-zt^{k_{2}}\right)^{-\alpha_{1}}\Theta\left(\kappa_{l};-\frac{b}{t}-\frac{d}{1-t}\right)dt.
\end{multline}
\[
\Re\left(\beta_{1}\right)>\Re\left(\alpha_{2}\right)>0; \Re\left(b\right),\Re\left(d\right)>0; b=d=0, \left|\arg\left(1-z\right)<\pi\right|
\]

As an example, let us compute the case for $k_{1}=2$ and $z=1$,
\begin{align*}
 & \frac{\Gamma\left(\beta_{1}\right)}{\Gamma\left(\alpha_{2}\right)\Gamma\left(\beta_{1}-\alpha_{2}\right)}\int_{0}^{1}t^{\alpha_{2}-1}\left(1-t\right)^{\beta_{1}-\alpha_{2}-1}\left(1-t^{2}\right)^{-\alpha_{1}}\Theta\left(\kappa_{l};-\frac{b}{t}-\frac{d}{1-t}\right)dt\\
 & =\frac{\Gamma\left(\beta_{1}\right)}{\Gamma\left(\alpha_{2}\right)\Gamma\left(\beta_{1}-\alpha_{2}\right)}\int_{0}^{1}t^{\alpha_{2}-1}\left(1-t\right)^{\beta_{1}-\alpha_{2}-1}\left(1-t\right)^{-\alpha_{1}}\left(1+t\right)^{-\alpha_{1}}\Theta\left(\kappa_{l};-\frac{b}{t}-\frac{d}{1-t}\right)dt\\
 & =\frac{\Gamma\left(\beta_{1}\right)}{\Gamma\left(\alpha_{2}\right)\Gamma\left(\beta_{1}-\alpha_{2}\right)}\sum_{k=0}^{\infty}\dbinom{-\alpha_{1}}{k}\int_{0}^{1}t^{\alpha_{2}+k-1}\left(1-t\right)^{\beta_{1}-\alpha_{2}-\alpha_{1}-1}\Theta\left(\kappa_{l};-\frac{b}{t}-\frac{d}{1-t}\right)dt\\
 & =\frac{\Gamma\left(\beta_{1}\right)}{\Gamma\left(\alpha_{2}\right)\Gamma\left(\beta_{1}-\alpha_{2}\right)}\sum_{k=0}^{\infty}\frac{\left(-1\right)^{k}\left(\alpha_{1}\right)_{k}}{k!}\mathcal{B}_{b,d}^{\left(\kappa_{l}\right)}\left(\alpha_{2}+k,\beta_{1}-\alpha_{2}-\alpha_{1}\right)\\
 & =\frac{\Gamma\left(\beta_{1}\right)}{\Gamma\left(\alpha_{2}\right)\Gamma\left(\beta_{1}-\alpha_{2}\right)}\frac{\Gamma\left(\beta_{1}-\alpha_{2}-\alpha_{1}\right)\Gamma\left(\alpha_{2}\right)}{\Gamma\left(\beta_{1}-\alpha_{1}\right)}\sum_{k=0}^{\infty}\left(\alpha_{1}\right)_{k}\frac{\mathcal{B}_{b,d}^{\left(\kappa_{l}\right)}\left(\alpha_{2}+k,\beta_{1}-\alpha_{2}-\alpha_{1}\right)}{B\left(\alpha_{2},\beta_{1}-\alpha_{2}-\alpha_{1}\right)}\frac{\left(-1\right)^{k}}{k!}\\
 & =\frac{\Gamma\left(\beta_{1}\right)\Gamma\left(\beta_{1}-\alpha_{2}-\alpha_{1}\right)}{\Gamma\left(\beta_{1}-\alpha_{2}\right)\Gamma\left(\beta_{1}-\alpha_{1}\right)}{}_{2}F_{1}^{\left(\kappa_{l}\right)}\left[\begin{array}{c}
\alpha_{1},\alpha_{2}\\
\beta_{1}-\alpha_{1}
\end{array};-1;b,d\right].
\end{align*}
We conclude our result as the following theorem:
\begin{Theorem}
If $\Re\left(\beta_{1}\right)>\Re\left(\alpha_{2}\right)>0$
and $\Re\left(\beta_{1}-\alpha_{2}-\alpha_{1}\right)>0$, then
\begin{equation}
_{2}F_{1}^{\left(\kappa_{l}\right)}\left[\begin{array}{c}
\left(\alpha_{1},1\right),\left(\alpha_{2},2\right)\\
\beta_{1}
\end{array};1;b,d\right]=\frac{\Gamma\left(\beta_{1}\right)\Gamma\left(\beta_{1}-\alpha_{2}-\alpha_{1}\right)}{\Gamma\left(\beta_{1}-\alpha_{2}\right)\Gamma\left(\beta_{1}-\alpha_{1}\right)}{}_{2}F_{1}^{\left(\kappa_{l}\right)}\left[\begin{array}{c}
\alpha_{1},\alpha_{2}\\
\beta_{1}-\alpha_{1}
\end{array};-1;b,d\right].
\end{equation}
\end{Theorem}

\textbf{Theorem 2.8}  is in fact a generalization of the following well-known result \cite[p.117, Theorem 2.3]{An Integral Representation of Some Hypergeometric Functions}:
\begin{equation}
_{3}F_{2}\left[\begin{array}{c}
a,\frac{b}{2},\frac{b+1}{2}\\
\frac{c}{2},\frac{c+1}{2}
\end{array};1\right]=\frac{\Gamma\left(c\right)\Gamma\left(c-a-b\right)}{\Gamma\left(c-b\right)\Gamma\left(c-a\right)}{}_{2}F_{1}\left[\begin{array}{c}
a,b\\
c-a
\end{array};-1\right].
\end{equation}
It also show some intrinsic relationships between our extended hypergeometric
function and the ordinary generalized hypergeometric functions.

In general, we have for definition (37) the following integral representation, which is a generalization of \cite[p.104, Theorem 38]{Special Functions Rainville}.
\begin{Theorem}
If $p\leq q+1$, $\Re\left(\beta_{q}\right)>\Re\left(\alpha_{p}\right)>0$, $k_j\in\mathbb{Z}^{+}\cup\{0\}, j=1,\cdots,p$ and none of $\beta_{1},\dots,\beta_{q}$ is zero or a negative integer, we have (inside the region of convergence)
\begin{align}
_{p}F_{q}^{\left(\kappa_{l}\right)}\left[\begin{array}{ccc}
\left(\alpha_{1},k_{1}\right) & \cdots & \left(\alpha_{p},k_{p}\right)\\
\beta_{1} & \cdots & \beta_{q}
\end{array};cz^{k_{p}};b,d\right] & =\frac{\Gamma\left(\beta_{q}\right)z^{1-\beta_{q}}}{\Gamma\left(\alpha_{p}\right)\Gamma\left(\beta_{q}-\alpha_{p}\right)}\int_{0}^{z}t^{\alpha_{p}-1}\left(z-t\right)^{\beta_{q}-a_{p}-1}\nonumber \\
 & \ \ \times{}_{p-1}F_{q-1}^{\left(\kappa_{l}\right)}\left[\begin{array}{ccc}
\left(\alpha_{1},k_{1}\right) & \cdots & \left(\alpha_{p-1},k_{p-1}\right)\\
\beta_{1} & \cdots & \beta_{q-1}
\end{array};ct^{k_{p}};b,d\right]\nonumber \\
 & \ \ \times\Theta\left(\kappa_{l};-\frac{bz}{t}-\frac{dz}{z-t}\right)dt.
\end{align}
\end{Theorem}
\begin{Proof}Denote the right-hand side of (42) by $A\left(z\right)$. As indicated before, we just need to prove the case that $p=q+1$. Let $t=zv$ in $A\left(z\right)$. Then we have
\begin{align*}
A\left(z\right) & =\frac{\Gamma\left(\beta_{q}\right)z^{\beta_{q}-1}}{\Gamma\left(\alpha_{p}\right)\Gamma\left(\beta_{q}-\alpha_{p}\right)}\int_{0}^{1}v^{\alpha_{q+1}-1}\left(1-v\right)^{\beta_{q}-a_{q+1}-1}\\
 & \ \ \times{}_{q}F_{q-1}^{\left(\kappa_{l}\right)}\left[\begin{array}{ccc}
\left(\alpha_{1},k_{1}\right) & \cdots & \left(\alpha_{q},k_{q}\right)\\
\beta_{1} & \cdots & \beta_{q-1}
\end{array};c\left(zv\right)^{k_{q+1}};b,d\right]\Theta\left(\kappa_{l};-\frac{b}{v}-\frac{d}{1-v}\right)dv\\
 & =\frac{\Gamma\left(\beta_{q}\right)z^{\beta_{q}-1}}{\Gamma\left(\alpha_{q+1}\right)\Gamma\left(\beta_{q}-\alpha_{q+1}\right)}\sum_{m=0}^{\infty}\left(\alpha_{1}\right)_{k_{1}m}\prod_{j=1}^{q-1}\frac{\mathcal{B}_{b,d}^{\left(\kappa_{l}\right)}\left(\alpha_{j+1}+k_{j+1}m,\beta_{j}-\alpha_{j+1}\right)}{B\left(\alpha_{j+1},\beta_{j}-\alpha_{j+1}\right)}\frac{\left(cz^{k_{q+1}}\right)^{m}}{m!}\\
 & \ \ \times\int_{0}^{1}v^{\alpha_{q+1}+k_{q+1}m-1}\left(1-v\right)^{\beta_{q}-a_{q+1}-1}\Theta\left(\kappa_{l};-\frac{b}{v}-\frac{d}{1-v}\right)dv\\
 & =z^{\beta_{q}-1}\sum_{m=0}^{\infty}\left(\alpha_{1}\right)_{k_{1}m}\prod_{j=1}^{q}\frac{\mathcal{B}_{b,d}^{\left(\kappa_{l}\right)}\left(\alpha_{j+1}+k_{j+1}m,\beta_{j}-\alpha_{j+1}\right)}{B\left(\alpha_{j+1},\beta_{j}-\alpha_{j+1}\right)}\frac{\left(cz^{k_{q+1}}\right)^{m}}{m!}\\
 & =z^{\beta_{q}-1}{}_{q+1}F_{q}^{\left(\kappa_{l}\right)}\left[\begin{array}{cccc}
\left(\alpha_{1},k_{1}\right) & \cdots & \left(\alpha_{q},k_{q}\right) & \left(\alpha_{q+1},k_{q+1}\right)\\
\beta_{1} & \cdots & \beta_{q-1} & \beta_{q}
\end{array};cz^{k_{q+1}};b,d\right].
\end{align*}
\end{Proof}

\textbf{Theorem 2.9} embody the significance and advantages of our modified definition (37) of extended generalized hypergeometric function. It not only enables us to rebuild a lot of classical results of hypergeometric function in a new way but also provide some new important interpretations. Here is one of them. In \cite{A Class of Extended Fractional Derivative Operators
and Associated Generating Relations Involving Hypergeometric Functions}, H. M. Srivastava definite the following fractional derivative operator:
\begin{equation}
\mathcal{D}_{z,\left(\kappa_{l}\right)}^{\mu,b,d}\left\{ f\left(z\right)\right\} =\begin{cases}
\frac{1}{\Gamma\left(-\mu\right)}\int_{0}^{z}\left(z-t\right)^{-\mu-1}\Theta\left(\kappa_{l};-\frac{bz}{t}-\frac{dz}{z-t}\right)f\left(t\right)dt & \left(\Re\left(\mu\right)<0\right)\\[15pt]
\frac{d^{m}}{dz^{m}}\left\{ \mathcal{D}_{z,\left(\kappa_{l}\right)}^{\mu-m,b,d}\left\{ f\left(z\right)\right\} \right\}  & \left(m-1\leq\Re\left(\mu\right)\leq m\left(m\in\mathbb{N}\right)\right)
\end{cases}
\end{equation}
The path of integration in (43) is a line in the complex $t-$plane from $t=0$ to $t=z$.

By using operator (43) we can rewrite (42) as
\begin{align}
 & _{p}F_{q}^{\left(\kappa_{l}\right)}\left[\begin{array}{ccc}
\left(\alpha_{1},k_{1}\right) & \cdots & \left(\alpha_{p},k_{p}\right)\\
\beta_{1} & \cdots & \beta_{q}
\end{array};cz^{k_{p}};b,d\right]\nonumber \\
 & =\frac{\Gamma\left(\beta_{q}\right)}{\Gamma\left(\alpha_{p}\right)}z^{1-\beta_{q}}\mathcal{D}_{z,\left(\kappa_{l}\right)}^{-\left(\beta_{q}-\alpha_{p}\right),b,d}\left\{ t^{\alpha_{p}-1}{}_{p-1}F_{q-1}^{\left(\kappa_{l}\right)}\left[\begin{array}{ccc}
\left(\alpha_{1},k_{1}\right) & \cdots & \left(\alpha_{p-1},k_{p-1}\right)\\
\beta_{1} & \cdots & \beta_{q-1}
\end{array};ct^{k_{p}};b,d\right]\right\} .
\end{align}

\section{{\normalsize Extended Appell and Lauricella Hypergeometric Functions}}
In this section, we focus on some extended multivariate hypergeometric functions, i.e., extended Appell's hypergeometric functions and extended Lauricella's hypergeometric functions.

For the sake of clarity and easy readability, we may first study the properties of extended Appell's hypergeometric functions and then
we can view the extended Lauricella's hypergeometric functions as a further generalization of Appell functions.
\begin{Definition}
\cite{A Class of Extended Fractional Derivative Operators
and Associated Generating Relations Involving Hypergeometric Functions}: The extended Appell's hypergeometric functions
\[
\mathbf{F}_{1}^{\left(\kappa_{l}\right)}\left(\alpha,\beta_{1},\beta_{2};\gamma_{1};x,y;b,d\right)\ and\ \mathbf{F}_{2}^{\left(\kappa_{l}\right)}\left(\alpha,\beta_{1},\beta_{2};\gamma_{1},\gamma_{2};x,y;b,d\right)
\]
in two variables are defined by
\begin{equation}
\mathbf{F}_{1}^{\left(\kappa_{l}\right)}\left(\alpha,\beta_{1},\beta_{2};\gamma_{1};x,y;b,d\right)=\sum_{m,n=0}^{\infty}\left(\beta_{1}\right)_{m}\left(\beta_{2}\right)_{n}\frac{\mathcal{B}_{b,d}^{\left(\kappa_{l}\right)}\left(\alpha+m+n,\gamma_{1}-\alpha\right)}{B\left(\alpha,\gamma_{1}-\alpha\right)}\frac{x^{m}y^{n}}{m!n!},
\end{equation}
\[
\left(\max\left\{ \left|x\right|,\left|y\right|\right\} <1;\min\left\{ \Re\left(b\right),\Re\left(d\right)\right\} \geq0\right)
\]
\begin{equation}
\mathbf{F}_{2}^{\left(\kappa_{l}\right)}\left(\alpha,\beta_{1},\beta_{2};\gamma_{1},\gamma_{2};x,y;b,d\right)=\sum_{m,n=0}^{\infty}\left(\alpha\right)_{m+n}\frac{\mathcal{B}_{b,d}^{\left(\kappa_{l}\right)}\left(\beta_{1}+m,\gamma_{1}-\beta_{1}\right)B_{b,d}^{\left(\kappa_{l}\right)}\left(\beta_{2}+n,\gamma_{2}-\beta_{2}\right)}{B\left(\beta_{1},\gamma_{1}-\beta_{1}\right)B\left(\beta_{2},\gamma_{2}-\beta_{2}\right)}\frac{x^{m}y^{n}}{m!n!}.
\end{equation}
\[
\left(\left|x\right|+\left|y\right|<1;\min\left\{ \Re\left(b\right),\Re\left(d\right)\right\} \geq0\right)
\]
\end{Definition}

The following Euler type integral representations for Extended Appell hypergeometric functions $\mathbf{F}_{1}^{\left(\kappa_{l}\right)}$ and $\mathbf{F}_{2}^{\left(\kappa_{l}\right)}$ are easily obtained and their proofs will be omitted.
\begin{Theorem}
For the extended Appell functions defined by (45) and (46), the following integral representations hold true:
\begin{multline}
\mathbf{F}_{1}^{\left(\kappa_{l}\right)}\left(\alpha,\beta_{1},\beta_{2};\gamma_{1};x,y;b,d\right)\\
=\frac{1}{B\left(\alpha,\gamma-\alpha\right)}
\int_{0}^{1}t^{\alpha-1}\left(1-t\right)^{\gamma_{1}-\alpha-1}\left(1-xt\right)^{-\beta_{1}}\left(1-yt\right)^{-\beta_{2}}\Theta\left(\kappa_{l};-\frac{b}{t}-\frac{d}{1-t}\right)dt,
\end{multline}
\[
\left(\min\left\{ \Re\left(b\right),\Re\left(d\right)\right\} \geq0;\max\left\{ \left|\arg\left(1-x\right)\right|,\left|\arg\left(1-y\right)\right|\right\} <\pi;\Re\left(\gamma\right)>\Re\left(\alpha\right)>0\right)
\]
\begin{multline}
\mathbf{F}_{2}^{\left(\kappa_{l}\right)}\left(\alpha,\beta_{1},\beta_{2};\gamma_{1},\gamma_{2};x,y;b,d\right)=\frac{1}{B\left(\beta_{1},\gamma_{1}-\beta_{1}\right)}\frac{1}{B\left(\beta_{2},\gamma_{2}-\beta_{2}\right)}\\
\int_{0}^{1}\int_{0}^{1}\frac{t^{\beta_{1}-1}\left(1-t\right)^{\gamma_{1}-\beta_{1}-1}s^{\beta_{2}-1}\left(1-s\right)^{\gamma_{2}-\beta_{2}-1}}{\left(1-xt-ys\right)^{\alpha}}\Theta\left(\kappa_{l};-\frac{b}{t}-\frac{d}{1-t}\right)\Theta\left(\kappa_{l};-\frac{b}{s}-\frac{d}{1-s}\right)dtds.
\end{multline}
\[
\left(\min\left\{ \Re\left(b\right),\Re\left(d\right)\right\} \geq0;\left|x\right|+\left|y\right|<1;\Re\left(\gamma_{1}\right)>\Re\left(\beta_{1}\right)>0,\Re\left(\gamma_{2}\right)>\Re\left(\beta_{2}\right)>0\right)
\]
\end{Theorem}

By using the integral representations of $\mathbf{F}_{1}^{\left(\kappa_{l}\right)}\left(\alpha,\beta,\beta';\gamma;x,y;b,d\right)$
and $\mathbf{F}_{2}^{\left(\kappa_{l}\right)}\left(\alpha,\beta,\beta';\gamma,\gamma';x,y;b,d\right)$, we can derive some transformation formulas for them.
\begin{Theorem}
For the extended Appell hypergeometric function $\mathbf{F}_{1}^{\left(\kappa_{l}\right)}$, we have the following transformation formula:
\begin{equation}
\mathbf{F}_{1}^{\left(\kappa_{l}\right)}\left(\alpha,\beta,\beta';\gamma;x,y;b,d\right)=\left(1-x\right)^{-\beta}\left(1-y\right)^{-\beta'}\mathbf{F}_{1}^{\left(\kappa_{l}\right)}\left(\alpha,\beta,\beta';\gamma;\frac{x}{x-1},\frac{y}{y-1};d,b\right).
\end{equation}
\end{Theorem}
\begin{Proof} Let us put $t=1-s$. Then
\begin{align*}
\mathbf{F}_{1}^{\left(\kappa_{l}\right)}\left(\alpha,\beta,\beta';\gamma;x,y;b,d\right) & =\frac{1}{B\left(\alpha,\gamma-\alpha\right)}\int_{0}^{1}s^{\gamma-\alpha-1}\left(1-s\right)^{\alpha-1}\left(1-x\left(1-s\right)\right)^{-\beta}\left(1-y\left(1-s\right)\right)^{-\beta'}\\
 & \ \ \times\Theta\left(\kappa_{l};-\frac{b}{1-s}-\frac{d}{s}\right)ds\\
 & =\frac{1}{B\left(\alpha,\gamma-\alpha\right)}\int_{0}^{1}s^{\gamma-\alpha-1}\left(1-s\right)^{\alpha-1}\left(1-x+xs\right)^{-\beta}\left(1-y+ys\right)^{-\beta'}\\
 & \ \ \times\Theta\left(\kappa_{l};-\frac{b}{1-s}-\frac{d}{s}\right)ds\\
 & =\frac{\left(1-x\right)^{-\beta}\left(1-y\right)^{\beta'}}{B\left(\alpha,\gamma-\alpha\right)}\int_{0}^{1}s^{\gamma-\alpha-1}\left(1-s\right)^{\alpha-1}\left(1-\frac{x}{x-1}s\right)^{-\beta}\left(1-\frac{y}{y-1}s\right)^{-\beta'}\\
 & \times\Theta\left(\kappa_{l};-\frac{d}{s}-\frac{b}{1-s}\right)ds
\end{align*}
such that
\[
\mathbf{F}_{1}^{\left(\kappa_{l}\right)}\left(\alpha,\beta,\beta';\gamma;x,y;b,d\right)=\left(1-x\right)^{-\beta}\left(1-y\right)^{\beta'}\mathbf{F}_{1}^{\left(\kappa_{l}\right)}\left(\gamma-\alpha,\beta,\beta';\gamma;\frac{x}{x-1},\frac{y}{y-1};d,b\right).
\]
\end{Proof}
\begin{Remark}
The classical theory of Appell hypergeometric functions have
found that integral
\[
\int_{0}^{1}u^{\alpha-1}\left(1-u\right)^{\gamma-\alpha-1}\left(1-ux\right)^{-\beta}\left(1-uy\right)^{-\beta'}du
\]
keeps unchanged under the following five transformations \cite[p.195]{SPECIAL FUNCTIONS}
\[
\left.\begin{array}{cc}
u=1-v,\\
u=\frac{v}{1-y+vy}, & u=\frac{1-v}{1-vx},\\
u=\frac{v}{1-x+vx}, & y=\frac{1-v}{1-vy}.
\end{array}\right\}
\]
However, only the first transformation $u=1-v$ enable us to get a simple extension.
\end{Remark}
\begin{Theorem} For $b=d$ we have the following transformation formulas for Appell hypergeometric function $\mathbf{F}_{2}^{\left(\kappa_{l}\right)}$:
\begin{equation}
\mathbf{F}_{2}^{\left(\kappa_{l}\right)}\left(\alpha,\beta,\beta';\gamma,\gamma';x,y;b\right)=\left(1-x\right)^{-\alpha}\mathbf{F}_{2}^{\left(\kappa_{l}\right)}\left(\alpha,\gamma-\beta,\beta';\gamma,\gamma';\frac{x}{x-1},\frac{-y}{x-1};b\right),
\end{equation}
\begin{equation}
\mathbf{F}_{2}^{\left(\kappa_{l}\right)}\left(\alpha,\beta,\beta';\gamma,\gamma';x,y;b\right)=\left(1-y\right)^{-\alpha}\mathbf{F}_{2}^{\left(\kappa_{l}\right)}\left(\alpha,\beta,\gamma'-\beta';\gamma,\gamma';\frac{-x}{y-1},\frac{y}{y-1};b\right),
\end{equation}
\begin{equation}
\mathbf{F}_{2}^{\left(\kappa_{l}\right)}\left(\alpha,\beta,\beta';\gamma,\gamma';x,y;b\right)=\left(1-x-y\right)^{-\alpha}\mathbf{F}_{2}^{\left(\kappa_{l}\right)}\left(\alpha,\gamma-\beta,\gamma'-\beta';\gamma,\gamma';\frac{-x}{1-x-y},\frac{-y}{1-x-y};b\right).
\end{equation}
\end{Theorem}
\begin{Proof}
By using the double integral representation of \emph{$\mathbf{F}_{2}^{\left(\kappa_{l}\right)}$}
with $b=d$ and the transformations
\[
\left.\begin{array}{ccc}
(a) & t=1-t', & s=s';\\
(b) & t=t' & s=s-1;\\
(c) & t=1-t', & s=1-s';
\end{array}\right\}
\]
we can directly obtain the transformation formulas for \emph{$\mathbf{F}_{2}^{\left(\kappa_{l}\right)}\left(\alpha,\beta,\beta';\gamma,\gamma';x,y;b\right)$},
respectively.

Note that formula (52) is also valid for more general case, namely,
\begin{equation}
\mathbf{F}_{2}^{\left(\kappa_{l}\right)}\left(\alpha,\beta,\beta';\gamma,\gamma';x,y;b,d\right)=\left(1-x-y\right)^{-\alpha}\mathbf{F}_{2}^{\left(\kappa_{l}\right)}\left(\alpha,\gamma-\beta,\gamma'-\beta';\gamma,\gamma';\frac{-x}{1-x-y},\frac{-y}{1-x-y};d,b\right).
\end{equation}
\end{Proof}

\subsection{{\normalsize \textbf{Other Integral Representations of Appell's Hypergeometric
Functions} $\mathbf{F}_{1}^{\left(\kappa_{l}\right)}$ \textbf{and} $\mathbf{F}_{2}^{\left(\kappa_{l}\right)}$}}
We now start to establish the Mellin-Barnes type integral for Appell functions $\mathbf{F}_{1}^{\left(\kappa_{l}\right)}$ and $\mathbf{F}_{2}^{\left(\kappa_{l}\right)}$. The tool we apply here is called \emph{Method of Bracket} \cite{Ramanujan Master Theorem}, which can be view as a multidimensional extension to Ramanujan's Master Theorem. This method enable us to evaluate multidimensional Mellin transform directly and efficiently. In this part, we will prove the following theorem with this method. Here we will put focus on the proving instead of explaining the method in detail. Because the method is easy-understanding, we will give a brief introduction of it when proving \textbf{Theorem 3.4.} More detailed information may be found in \cite{Ramanujan Master Theorem} and \cite{Definite integrals by the method of brackets}.
\begin{Theorem} For suitable parameters, the Mellin-Barnes integral representation of $\mathbf{F}_{1}^{\left(\kappa_{l}\right)}$ is
\begin{align}
\mathbf{F}_{1}^{\left(\kappa_{l}\right)}\left(\alpha,\beta_{1},\beta_{2};\gamma_{1};x,y;b,d\right) & =\frac{1}{\left(2\pi i\right)^{2}}\frac{\Gamma\left(\gamma_{1}\right)}{\Gamma\left(\beta_{1}\right)\Gamma\left(\beta_{2}\right)\Gamma\left(\alpha\right)\Gamma\left(\gamma_{1}-\alpha\right)}\int_{-i\infty}^{+i\infty}\int_{-i\infty}^{+i\infty}\mathcal{B}_{b,d}^{\left(\kappa_{l}\right)}\left(\alpha-s_{1}-s_{2},\gamma_{1}-\alpha\right)\nonumber \\
 & \ \ \times\Gamma\left(\beta_{1}-s_{1}\right)\Gamma\left(\beta_{2}-s_{2}\right)\Gamma\left(s_{1}\right)\Gamma\left(s_{2}\right)\left(-x\right)^{-s_{1}}\left(-y\right)^{-s_{2}}ds_{1}ds_{2}.
\end{align}
\end{Theorem}
\begin{Proof}
We first take the Mellin transform of $\mathbf{F}_{1}^{\left(\kappa_{l}\right)}\left(\alpha,\beta_{1},\beta_{2};\gamma_{1};-x,-y;b,d\right)$,
i.e.,
\begin{multline}
\int_{0}^{\infty}\int_{0}^{\infty}x^{s_{1}-1}y^{s_{2}-1}\mathbf{F}_{1}^{\left(\kappa_{l}\right)}\left(\alpha,\beta_{1},\beta_{2};\gamma_{1};-x,-y;b,d\right)dxdy =\sum_{m,n=0}^{\infty}\left(\beta_{1}\right)_{m}\left(\beta_{2}\right)_{n}\frac{\mathcal{B}_{b,d}^{\left(\kappa_{l}\right)}\left(\alpha+m+n,\gamma_{1}-\alpha\right)}{B\left(\alpha,\gamma_{1}-\alpha\right)}\\
\times\frac{\left(-1\right)^{m}}{\Gamma\left(m+1\right)}\frac{\left(-1\right)^{n}}{\Gamma\left(n+1\right)}\int_{0}^{\infty}x^{s_{1}+m-1}dx\int_{0}^{\infty}y^{s_{2}+n-1}dy
=\sum_{m,n=0}^{\infty}\phi_{m}\phi_{n}f\left(m,n\right)\left\langle s_{1}+m\right\rangle \left\langle s_{2}+n\right\rangle
\end{multline}
where
\[
f\left(m,n\right)=\left(\beta_{1}\right)_{m}\left(\beta_{2}\right)_{n}\frac{\mathcal{B}_{b,d}^{\left(\kappa_{l}\right)}\left(\alpha+m+n,\gamma_{1}-\alpha\right)}{B\left(\alpha,\gamma_{1}-\alpha\right)}=\frac{\Gamma\left(\beta_{1}+m\right)}{\Gamma\left(\beta_{1}\right)}\frac{\Gamma\left(\beta_{2}+n\right)}{\Gamma\left(\beta_{2}\right)}\frac{\mathcal{B}_{b,d}^{\left(\kappa_{l}\right)}\left(\alpha+m+n,\gamma_{1}-\alpha\right)}{B\left(\alpha,\gamma_{1}-\alpha\right)},
\]
\[
\phi_m=\frac{(-1)^m}{\Gamma(m+1)}, \phi_n=\frac{(-1)^n}{\Gamma(n+1)},
\]
\[
\left\langle s_{1}+m\right\rangle=\int_{0}^{\infty}x^{s_{1}+m-1}dx,\  \left\langle s_{2}+n\right\rangle=\int_{0}^{\infty}y^{s_{2}+n-1}dy.
\]
The symbol $\phi_m$ and $\phi_n$ are called the \emph{indicator} of $m$ and $n$, respectively \cite[p.53, Definition 2.2]{Definite integrals by the method of brackets}. Then use the \emph{method of brackets} \textbf{Rule III} \cite{Ramanujan Master Theorem} we have
\begin{multline}
\int_{0}^{\infty}\int_{0}^{\infty}x^{s_{1}-1}y^{s_{2}-1}\mathbf{F}_{1}^{\left(\kappa_{l}\right)}\left(\alpha,\beta_{1},\beta_{2};\gamma_{1};-x,-y;b,d\right)dxdy =\frac{1}{\left|\det\left(A\right)\right|}\Gamma\left(-m^{*}\right)\Gamma\left(-n^{*}\right)f\left(m^{*},n^{*}\right)\\
=\frac{\Gamma\left(\beta_{1}-s_{1}\right)}{\Gamma\left(\beta_{1}\right)}\frac{\Gamma\left(\beta_{2}-s_{2}\right)}{\Gamma\left(\beta_{2}\right)}\frac{\mathcal{B}_{b,d}^{\left(\kappa_{l}\right)}\left(\alpha-s_{1}-s_{2},\gamma_{1}-\alpha\right)}{B\left(\alpha,\gamma_{1}-\alpha\right)}\Gamma\left(s_{1}\right)\Gamma\left(s_{2}\right),
\end{multline}
where $(m^{*},n^{*})$ is obtained by vanishing the expressions in the brackets $\left\langle s_{1}+m\right\rangle$, $\left\langle s_{2}+n\right\rangle$. In general, $(m^{*},n^{*})$ is the unique solution to the linear system obtained by vanishing the expressions in the brackets. And the determinant of its coefficient matrix is denoted by $\det(A)$. Here, we can see $\left\langle a\right\rangle$ is only a symbol associated with divergent integrals  $\int_0^{\infty}x^{a-1}dx$. But by introducing such a symbol, we in fact reduce the evaluation of definite integral to solving a linear system of equations, which is the basic idea of the \emph{method of brackets}.

Now we can take the inverse Mellin transform of (56) to get
\begin{align*}
\mathbf{F}_{1}^{\left(\kappa_{l}\right)}\left(\alpha,\beta_{1},\beta_{2};\gamma_{1};-x,-y;b,d\right) & =\frac{1}{\left(2\pi i\right)^{2}}\int_{-i\infty}^{+i\infty}\int_{-i\infty}^{+i\infty}\frac{\Gamma\left(\beta_{1}-s_{1}\right)}{\Gamma\left(\beta_{1}\right)}\frac{\Gamma\left(\beta_{2}-s_{2}\right)}{\Gamma\left(\beta_{2}\right)}\frac{\mathcal{B}_{b,d}^{\left(\kappa_{l}\right)}\left(\alpha-s_{1}-s_{2},\gamma_{1}-\alpha\right)}{B\left(\alpha,\gamma_{1}-\alpha\right)}\\
 & \ \ \times\Gamma\left(s_{1}\right)\Gamma\left(s_{2}\right)x^{-s_{1}}y^{-s_{2}}ds_{1}ds_{2}.
\end{align*}
Theorem 3.4 follows from replacing $-x$ and $-y$ with $x$ and $y$.
\end{Proof}

The Mellin-barnes double contour integral representation for extended
Appell's hypergeometric function $\mathbf{F}_{2}^{\left(\kappa_{l}\right)}\left(\alpha,\beta_{1},\beta_{2};\gamma_{1},\gamma_{2};x,y;b,d\right)$
is given in the following theorem. We prove this just the same way when proving \textbf{Theorem 3.4}.
\begin{Theorem}
For suitable parameters, the Mellin-Barnes integral representation of $\mathbf{F}_{2}^{\left(\kappa_{l}\right)}$ is
\begin{multline}
\mathbf{F}_{2}^{\left(\kappa_{l}\right)}\left(\alpha,\beta_{1},\beta_{2};\gamma_{1},\gamma_{2};x,y;b,d\right) =\frac{1}{\left(2\pi i\right)^{2}}\frac{\Gamma\left(\gamma_{1}\right)\Gamma\left(\gamma_{2}\right)}{\Gamma\left(\beta_{1}\right)\Gamma\left(\gamma_{1}-\beta_{1}\right)\Gamma\left(\beta_{2}\right)\Gamma\left(\gamma_{2}-\beta_{2}\right)\Gamma\left(\alpha\right)}\\
\int_{-i\infty}^{i\infty}\int_{-i\infty}^{i\infty}\mathcal{B}_{b,d}^{\left(\kappa_{l}\right)}\left(\beta_{1}-s_{1},\gamma_{1}-\beta_{1}\right)\nonumber B_{b,d}^{\left(\kappa_{l}\right)}\left(\beta_{2}-s_{2},\gamma_{2}-\beta_{2}\right)\Gamma\left(\alpha-s_{1}-s_{2}\right)\Gamma\left(s_{1}\right)\Gamma\left(s_{2}\right)\left(-x\right)^{-s_{1}}\left(-y\right)^{-s_{2}}ds_{1}ds_{2}.
\end{multline}
\end{Theorem}

\subsection{{\normalsize \textbf{Recursion formulas for Appell hypergeometric function}
$\mathbf{F}_{2}^{\left(\kappa_{l}\right)}$}}
Recursion formulas for classical Appell hypergeometric function $\mathbf{F}_{2}$
have been studied in \cite{Recursion formulas for Appells hypergeometric function with some applications to radiation field problems} by using contiguous function relations of the classical Gauss hypergeometric series $_{2}F_{1}$. These recursion
formulas can be used to evaluate the radiation field integrals. Now
we try to establish the similar formulas for our new Appell hypergeometric
function $\mathbf{F}_{2}^{\left(\kappa_{l}\right)}$.
\begin{Theorem} The extended Appell hypergeometric function $\mathbf{F}_{2}^{\left(\kappa_{l}\right)}$satisfies
the following identities:
\begin{multline}
\mathbf{F}_{2}^{\left(\kappa_{l}\right)}\left(\alpha,\beta_{1},\beta_{2}+n;\gamma_{1},\gamma_{2};x,y;b,d\right)\\
=\frac{\left(\gamma_{2}-\beta_{2}\right)_{2n}}{\left(\gamma_{2}-\beta_{2}\right)_{n}\left(\beta_{2}\right)_{n}}\sum_{i=1}^{n}\left(-n\right)_{i}\frac{\left(\beta_{2}\right)_{i+n}}{\left(\gamma_{2}\right)_{i+n}i!}\mathbf{F}_{2}^{\left(\kappa_{l}\right)}\left(\alpha,\beta_{1},\beta_{2}+n+i;\gamma_{1},\gamma_{2}+i+n;x,y;b,d\right)
\end{multline}
\[
\left(\left|x\right|+\left|y\right|<1;\min\left\{ \Re\left(b\right),\Re\left(d\right)\right\} \geq0;\ n\in\mathbb{N}_{0},\alpha,\beta_{1},\beta_{2}\in\mathbb{C};\gamma_{1},\gamma_{2}\in\mathbb{C}\backslash\mathbb{Z}_{0}^{-}\right)
\]
and
\begin{multline}
\mathbf{F}_{2}^{\left(\kappa_{l}\right)}\left(\alpha,\beta_{1},\beta_{2};\gamma_{1},\gamma_{2}+n;x,y;b,d\right)\\
=\frac{\left(\gamma_{2}\right)_{n}}{\left(\gamma_{2}-\beta_{2}\right)_{n}}\sum_{k=0}^{n}\left(-1\right)^{k}\left(\begin{array}{c}
n\\
k
\end{array}\right)\frac{\left(\beta_{2}\right)_{k}}{\left(\gamma_{2}\right)_{k}}\mathbf{F}_{2}^{\left(\kappa_{l}\right)}\left(\alpha,\beta_{1},\beta_{2}+k;\gamma_{1},\gamma_{2}+k;x,y;b,d\right).
\end{multline}
\[
\left(\left|x\right|+\left|y\right|<1;\min\left\{ \Re\left(b\right),\Re\left(d\right)\right\} \geq0;\ n\in\mathbb{N}_{0},\alpha,\beta_{1},\beta_{2}\in\mathbb{C};\gamma_{1},\gamma_{2}\in\mathbb{C}\backslash\mathbb{Z}_{0}^{-}\right)
\]
\end{Theorem}
\begin{Proof}
1. Using the integral representation (13) the following double integral representation
\begin{multline}
\mathbf{F}_{2}^{\left(\kappa_{l}\right)}\left(\alpha,\beta_{1},\beta_{2};\gamma_{1},\gamma_{2};x,y;b,d\right) =\frac{1}{B\left(\beta_{1},\gamma_{1}-\beta_{1}\right)}\frac{1}{B\left(\beta_{2},\gamma_{2}-\beta_{2}\right)}\\
\times\int_{0}^{1}\int_{0}^{1}\frac{t^{\beta_{1}-1}\left(1-t\right)^{\gamma_{1}-\beta_{1}-1}s^{\beta_{2}-1}\left(1-s\right)^{\gamma_{2}-\beta_{2}-1}}{\left(1-xt-ys\right)^{\alpha}}\nonumber \Theta\left(\kappa_{l};-\frac{b}{t}-\frac{d}{1-t}\right)\Theta\left(\kappa_{l};-\frac{b}{s}-\frac{d}{1-s}\right)dtds
\end{multline}
can be reduced to a single integral involving the extended Gauss hypergeometric function $_{2}F_{1}^{\left(\kappa_{l}\right)}$, i.e.,
\begin{multline}
\mathbf{F}_{2}^{\left(\kappa_{l}\right)}\left(\alpha,\beta_{1},\beta_{2};\gamma_{1},\gamma_{2};x,y;b,d\right) =\frac{1}{B\left(\beta_{1},\gamma_{1}-\beta_{1}\right)}\int_{0}^{1}t^{\beta_{1}-1}\left(1-t\right)^{\gamma_{1}-\beta_{1}-1}\left(1-xt\right)^{-\alpha}\\
{}_{2}F_{1}^{\left(\kappa_{l}\right)}\left(\begin{array}{c}
\alpha,\beta_{2}\\
\gamma_{2}
\end{array};\frac{y}{1-xt};b,d\right) \Theta\left(\kappa_{l};-\frac{b}{t}-\frac{d}{1-t}\right)dt.
\end{multline}
The assertion (57) then follows by substituting (30) into integral representation (59) with $\beta_{2}=\beta_{2}+n$. Thus
\begin{align}
& \mathbf{F}_{2}^{\left(\kappa_{l}\right)}\left(\alpha,\beta_{1},\beta_{2}+n;\gamma_{1},\gamma_{2};x,y;b,d\right)\nonumber \\
& =\frac{1}{B\left(\beta_{1},\gamma_{1}-\beta_{1}\right)}\int_{0}^{1}t^{\beta_{1}-1}\left(1-t\right)^{\gamma_{1}-\beta_{1}-1}\left(1-xt\right)^{-\alpha}{}_{2}F_{1}^{\left(\kappa_{l}\right)}\left(\begin{array}{c}
\alpha,\beta_{2}+n\\
\gamma_{2}
\end{array};\frac{y}{1-xt};b,d\right)\nonumber \\
 & \ \ \times\Theta\left(\kappa_{l};-\frac{b}{t}-\frac{d}{1-t}\right)dt\\
 & =\frac{1}{B\left(\beta_{1},\gamma_{1}-\beta_{1}\right)}\frac{\left(\gamma_{2}-\beta_{2}\right)_{2n}}{\left(\gamma_{2}-\beta_{2}\right)_{n}\left(\beta_{2}\right)_{n}}\sum_{i=1}^{n}\left(-n\right)_{i}\frac{\left(\beta_{2}\right)_{i+n}}{\left(\gamma_{2}\right)_{i+n}i!}\int_{0}^{1}t^{\beta_{1}-1}\left(1-t\right)^{\gamma_{1}-\beta_{1}-1}\left(1-xt\right)^{-\alpha}\nonumber \\
 & \ \ \times{}_{2}F_{1}^{\left(\kappa_{l}\right)}\left(\begin{array}{c}
\alpha,\beta_{2}+i+n\\
\gamma_{2}+i+n
\end{array};\frac{y}{1-xt};b,d\right)\Theta\left(\kappa_{l};-\frac{b}{t}-\frac{d}{1-t}\right)dt\\
 & =\frac{\left(\gamma_{2}-\beta_{2}\right)_{2n}}{\left(\gamma_{2}-\beta_{2}\right)_{n}\left(\beta_{2}\right)_{n}}\sum_{i=1}^{n}\left(-n\right)_{i}\frac{\left(\beta_{2}\right)_{i+n}}{\left(\gamma_{2}\right)_{i+n}i!}\mathbf{F}_{2}^{\left(\kappa_{l}\right)}\left(\alpha,\beta_{1},\beta_{2}+n+i;\gamma_{1},\gamma_{2}+i+n;x,y;b,d\right).
\end{align}

2. The proof of the second assertion of \textbf{Theorem 3.6} follows directly by substituting from the the assertion (29) of \textbf{Theorem 2.7} into the integral
representation (59) with $\gamma_2=\gamma_2+n$,
\begin{align}
& \mathbf{F}_{2}^{\left(\kappa_{l}\right)}\left(\alpha,\beta_{1},\beta_{2};\gamma_{1},\gamma_{2}+n;x,y;b,d\right)\nonumber \\& =\frac{1}{B\left(\beta_{1},\gamma_{1}-\beta_{1}\right)}\int_{0}^{1}t^{\beta_{1}-1}\left(1-t\right)^{\gamma_{1}-\beta_{1}-1}\left(1-xt\right)^{-\alpha}{}_{2}F_{1}^{\left(\kappa_{l}\right)}\left(\begin{array}{c}
\alpha,\beta_{2}\\
\gamma_{2}+n
\end{array};\frac{y}{1-xt};b,d\right)\nonumber \\
 & \ \ \times\Theta\left(\kappa_{l};-\frac{b}{t}-\frac{d}{1-t}\right)dt\nonumber \\
 & =\frac{\left(\gamma_{2}\right)_{n}}{\left(\gamma_{2}-\beta_{2}\right)_{n}}\sum_{k=0}^{n}\left(-1\right)^{k}\left(\begin{array}{c}
n\\
k
\end{array}\right)\frac{\left(\beta_{2}\right)_{k}}{\left(\gamma_{2}\right)_{k}}\int_{0}^{1}\frac{t^{\beta_{1}-1}\left(1-t\right)^{\gamma_{1}-\beta_{1}-1}}{B\left(\beta_{1},\gamma_{1}-\beta_{1}\right)}\left(1-xt\right)^{-\alpha}\nonumber \\
 & \ \ \times{}_{2}F_{1}^{\left(\kappa_{l}\right)}\left(\begin{array}{c}
\alpha,\beta_{2}+k\\
\gamma_{2}+k
\end{array};\frac{y}{1-xt};b,d\right)\Theta\left(\kappa_{l};-\frac{b}{t}-\frac{d}{1-t}\right)dt\\
 & =\frac{\left(\gamma_{2}\right)_{n}}{\left(\gamma_{2}-\beta_{2}\right)_{n}}\sum_{k=0}^{n}\left(-1\right)^{k}\left(\begin{array}{c}
n\\
k
\end{array}\right)\frac{\left(\beta_{2}\right)_{k}}{\left(\gamma_{2}\right)_{k}}\mathbf{F}_{2}^{\left(\kappa_{l}\right)}\left(\alpha,\beta_{1},\beta_{2}+k;\gamma_{1},\gamma_{2}+k;x,y;b,d\right).\nonumber
\end{align}
\end{Proof}
\begin{Remark}
The second assertion (58) of \textbf{Theorem 3.6} is clear a extension of
the formula given in \cite[p.549, Theorem 3]{Recursion formulas for Appells hypergeometric function with some applications to radiation field problems} and the first formula (57) seem to be new. In addition, the results obtained here are also
applicable to the extended Appell hypergeometric functions
\[
\mathbf{F}_{2}^{\left(\kappa_{l}\right)}\left(\alpha,\beta_{1},\beta_{2};\gamma_{1}+n,\gamma_{2};x,y;b,d\right)\ and\ \mathbf{F}_{2}^{\left(\kappa_{l}\right)}\left(\alpha,\beta_{1}+n,\beta_{2};\gamma_{1},\gamma_{2};x,y;b,d\right).
\]
It seems to be a little difficult to establish some recursion relations
for the functions:
\[
\mathbf{F}_{2}^{\left(\kappa_{l}\right)}\left(\alpha,\beta_{1},\beta_{2};\gamma_{1}-n,\gamma_{2};x,y;b,d\right)\ and\ \mathbf{F}_{2}^{\left(\kappa_{l}\right)}\left(\alpha,\beta_{1}-n,\beta_{2};\gamma_{1},\gamma_{2};x,y;b,d\right).
\]
\end{Remark}
\subsection{{\normalsize \textbf{Finite sum representation of the extended Appell's hypergeometric
function} $\mathbf{F}_{1}^{\left(\kappa_{l}\right)}$}}
Recently, Cuyt et al. \cite{A finite sum representation of the Appell series} obtained a finite algebraic sum representation of the Appell hypergeometric function $\mathbf{F}_{1}$ in the case when the parameters $\alpha,\beta_{1},\beta_{2}$ and $\gamma_{1}$are positive integers with $\gamma_{1}>\alpha$. Motivated by their works, we show that their results can be extended to obtain a finite sum representation of extended Appell's hypergeometric function $\mathbf{F}_{1}^{\left(\kappa_{l}\right)}$.

Before we state our results, we need the following important lemma.
\newtheorem{Lemma}{Lemma}
\begin{Lemma} \cite{A finite sum representation of the Appell series}: Let $\alpha=\frac{y}{y-x}$, $\beta=\frac{x}{x-y}$ and $X=\left(1-ux\right)^{-1}$, $Y=\left(1-uy\right)^{-1}$. Then for any positive integers $s$ and $t$,
\begin{equation}
X^{s}Y^{t}=\alpha^{s}\sum_{j=0}^{t-1}\dbinom{j+s-1}{s-1}\beta^{j}Y^{t-j}+\beta^{t}\sum_{k=0}^{s-1}\dbinom{k+t-1}{t-1}\alpha^{k}X^{s-k}.
\end{equation}
\end{Lemma}
\begin{Theorem} For any non-negative integers $s, t$ and $x\neq y$, we
have for $\left|x\right|<1,$ $\left|y\right|<1$
\begin{align}
\mathbf{F}_{1}^{\left(\kappa_{l}\right)}\left(1,s+1,t+1;2;x,y;b,d\right) & =y^{s+1}\sum_{j=0}^{t-1}\dbinom{j+s}{s}\frac{\left(-x\right)^{j}}{\left(y-x\right)^{j+s+1}}{}_{2}F_{1}^{\left(\kappa_{l}\right)}\left[\begin{array}{c}
\left(t-j+1\right),1\\
2
\end{array};y;b,d\right]\nonumber \\
 & \ \ +x^{t+1}\sum_{k=0}^{s-1}\dbinom{k+t}{t}\frac{\left(-y\right)^{k}}{\left(y-x\right)^{k+t+1}}{}_{2}F_{1}^{\left(\kappa_{l}\right)}\left[\begin{array}{c}
\left(s-k+1\right),1\\
2
\end{array};x;b,d\right]\nonumber \\
 & \ \ +\dbinom{t+s}{s}\frac{\left(-1\right)^{t}x^{t}y^{s}}{\left(y-x\right)^{t+s+1}}\left\{ y{}_{2}F_{1}^{\left(\kappa_{l}\right)}\left[\begin{array}{c}
1,1\\
2
\end{array};y;b,d\right]+x{}_{2}F_{1}^{\left(\kappa_{l}\right)}\left[\begin{array}{c}
1,1\\
2
\end{array};x;b,d\right]\right\}.
\end{align}
\end{Theorem}
\begin{Proof}
Put $\alpha=1$, $\beta_{1}=s+1$, $\beta_{2}=t+1$ and $\gamma_{2}=2$
in the integral representation (47), we obtain
\begin{equation}
\mathbf{F}_{1}^{\left(\kappa_{l}\right)}\left(1,s+1,t+1;2;x,y;b,d\right)=\int_{0}^{1}\left(1-xu\right)^{-\left(s+1\right)}\left(1-yu\right)^{-\left(t+1\right)}\Theta\left(\kappa_{l};-\frac{b}{u}-\frac{d}{1-u}\right)du.
\end{equation}
By using \textbf{lemma 1} we have
\begin{equation}
\left(1-xu\right)^{-\left(s+1\right)}\left(1-yu\right)^{-\left(t+1\right)}=\alpha^{s+1}\sum_{j=0}^{t}\dbinom{j+s}{s}\beta^{j}\left(1-uy\right)^{-\left(t+1-j\right)}+\beta^{t+1}\sum_{k=0}^{s}\dbinom{k+t}{t}\alpha^{k}\left(1-ux\right)^{-\left(s+1-k\right)}.
\end{equation}
Substituting (67) into (66) one get
\begin{align*}
\mathbf{F}_{1}^{\left(\kappa_{l}\right)}\left(1,s+1,t+1;2;x,y;b,d\right) & =\alpha^{s+1}\sum_{j=0}^{t}\dbinom{j+s}{s}\beta^{j}\int_{0}^{1}\left(1-uy\right)^{-\left(t+1-j\right)}\Theta\left(\kappa_{l};-\frac{b}{u}-\frac{d}{1-u}\right)du\\
 & \ \ +\beta^{t+1}\sum_{k=0}^{s}\dbinom{k+t}{t}\alpha^{k}\int_{0}^{1}\left(1-ux\right)^{-\left(s+1-k\right)}\Theta\left(\kappa_{l};-\frac{b}{u}-\frac{d}{1-u}\right)du\\
 & =\alpha^{s+1}\sum_{j=0}^{t}\dbinom{j+s}{s}\beta^{j}{}_{2}F_{1}^{\left(\kappa_{l}\right)}\left[\begin{array}{c}
\left(t-j+1\right),1\\
2
\end{array};y;b,d\right]\\
 & \ \ +\beta^{t+1}\sum_{k=0}^{s}\dbinom{k+t}{t}\alpha^{k}{}_{2}F_{1}^{\left(\kappa_{l}\right)}\left[\begin{array}{c}
\left(s-k+1\right),1\\
2
\end{array};x;b,d\right]\\
 & =y^{s+1}\sum_{j=0}^{t-1}\dbinom{j+s}{s}\frac{\left(-x\right)^{j}}{\left(y-x\right)^{j+s+1}}{}_{2}F_{1}^{\left(\kappa_{l}\right)}\left[\begin{array}{c}
\left(t-j+1\right),1\\
2
\end{array};y;b,d\right]\\
 & \ \ +x^{t+1}\sum_{k=0}^{s-1}\dbinom{k+t}{t}\frac{\left(-y\right)^{k}}{\left(y-x\right)^{k+t+1}}{}_{2}F_{1}^{\left(\kappa_{l}\right)}\left[\begin{array}{c}
\left(s-k+1\right),1\\
2
\end{array};x;b,d\right]\\
 & \ \ +\dbinom{t+s}{s}\frac{\left(-1\right)^{t}x^{t}y^{s}}{\left(y-x\right)^{t+s+1}}\left\{ y{}_{2}F_{1}^{\left(\kappa_{l}\right)}\left[\begin{array}{c}
1,1\\
2
\end{array};y;b,d\right]+x{}_{2}F_{1}^{\left(\kappa_{l}\right)}\left[\begin{array}{c}
1,1\\
2
\end{array};x;b,d\right]\right\},
\end{align*}
where we have applied the integral representation of extended Gauss hypergeometric function.
\end{Proof}
\begin{Remark}
If we set $b=d=0$ , then (65) is reduced to the following finite sum representation \cite[p.215, Theorem 2.1a]{A finite sum representation of the Appell series}
\begin{align}
\mathbf{F}_{1}\left(1,s+1,t+1;2;x,y\right) & =\frac{\left(-x\right)^{t}y^{s}}{\left(y-x\right)^{s+t+1}}\dbinom{s+t}{s}\left[\ln\left(1-x\right)-\ln\left(1-y\right)\right]\nonumber \\
 & \ \ -\sum_{j=0}^{t-1}\dbinom{j+s}{s}\frac{\left(-x\right)^{j}y^{s}\left[1-\left(1-y\right)^{j-t}\right]}{\left(y-x\right)^{s+j+1}\left(t-j\right)}\nonumber \\
 & \ \ -\sum_{k=0}^{s-1}\dbinom{k+t}{t}\frac{x^{t}\left(-y\right)^{k}\left[1-\left(1-x\right)^{k-s}\right]}{\left(x-y\right)^{t+k+1}\left(s-k\right)}.
\end{align}
\end{Remark}
\begin{Remark}
This method can be used to find more complicated finite sum representations of extended Appell's hypergeometric series $\mathbf{F}_{1}^{\left(\kappa_{l}\right)}$.
For instance, other finite sum representations given by \cite[p.215, Theorem 2.1 b,c]{A finite sum representation of the Appell series} can indeed be considered analogously in a simple and straightforward manner.
\end{Remark}
\ \\[5pt]

It marks the end of the discussion of extended Appell's hypergeometric functions. Next, we will put focus on the Lauricell's hypergeometric functions. However, not all four Lauricella's functions can be generalized to their new forms.here, we will give a known generalized result.
\begin{Definition}\cite{A Class of Extended Fractional Derivative Operators and Associated Generating Relations Involving Hypergeometric Functions} Lauricella hypergeometric function $\mathbf{F}_{D,\left(\kappa_{l}\right)}^{\left(r\right)}$ is defined by
\begin{multline}
\mathbf{F}_{D,\left(\kappa_{l}\right)}^{\left(r\right)}\left(\alpha,\beta_{1},\cdots,\beta_{r};\gamma;x_{1},\cdots,x_{r};b,d\right)\\
=\sum_{m_{1},\cdots,m_{r}=0}^{\infty}\left(b_{1}\right)_{m_{1}}\cdots\left(b_{r}\right)_{m_{r}}\frac{\mathcal{B}_{b,d}^{\left(\kappa_{l}\right)}\left(\alpha+m_{1}+\cdots+m_{r},\gamma-\alpha\right)}{B\left(\alpha,\gamma-\alpha\right)}\frac{x_{1}^{m_{1}}}{m_{1}!}\cdots\frac{x_{r}^{m_{r}}}{m_{r}!}.
\end{multline}
\[
\left(\max\left\{ \left|x_{1}\right|,\cdots,\left|x_{r}\right|\right\} <1;\min\left\{ \Re\left(b\right),\Re\left(d\right)\right\} \geq0\right)
\]
\end{Definition}

Its Euler type integral representation of Lauricella hypergeometric function $\mathbf{F}_{D,\left(\kappa_{l}\right)}^{\left(r\right)}$ is stated in the following theorem.
\begin{Theorem} For the extended Lauricella hypergeometric function $\mathbf{F}_{D,\left(\kappa_{l}\right)}^{\left(r\right)}$ defined by (69), the following integral integral representation hold true:
\begin{multline}
\mathbf{F}_{D,\left(\kappa_{l}\right)}^{\left(r\right)}\left(\alpha,\beta_{1},\cdots,\beta_{r};\gamma;x_{1},\cdots,x_{r};b,d\right)\\
=\frac{\Gamma\left(\gamma\right)}{\Gamma\left(\alpha\right)\Gamma\left(\gamma-\alpha\right)}\int_{0}^{1}t^{\alpha-1}\left(1-t\right)^{\gamma-\alpha-1}\prod_{j=1}^{r}\left(1-x_{j}t\right)^{-\beta_{j}}\Theta\left(\kappa_{l};-\frac{b}{t}-\frac{d}{1-t}\right)dt.
\end{multline}
\[
\left(\Re\left(b\right),\Re\left(d\right)>0;b=d=0, \max\left\{ \left|\arg\left(1-x_{1}\right)\right|,\cdots\left|\arg\left(1-x_{r}\right)\right|\right\} <\pi;\Re\left(\gamma\right)>\Re\left(\alpha\right)>0\right)
\]
\end{Theorem}

As a direct consequence of above theorem, we have the following summation formula
\[
\mathbf{F}_{D,\left(\kappa_{l}\right)}^{\left(r\right)}\left(\alpha,\beta_{1},\cdots,\beta_{r};\gamma;1,\cdots,1;b,d\right)=\frac{\Gamma\left(\gamma\right)}{\Gamma\left(\alpha\right)\Gamma\left(\gamma-\alpha\right)}\mathcal{B}_{b,d}^{\left(\kappa_{l}\right)}\left(\alpha,\gamma-\alpha-\beta_{1}-\cdots-\beta_{r}\right).
\]
If we put $x_{1}=x_{2}=\cdots=x_{r}=x$,
\[
\mathbf{F}_{D,\left(\kappa_{l}\right)}^{\left(r\right)}\left(\alpha,\beta_{1},\cdots,\beta_{r};\gamma;x,\cdots,x;b,d\right)={}_{2}F_{1}^{\left(\kappa_l\right)}\left(\alpha;\beta_{1}+\cdots+\beta_{r};\gamma;x;b,d\right).
\]
\begin{Theorem} One has
\begin{align}
\int_{a}^{b}\left(t-a\right)^{\alpha-1}\left(b-t\right)^{\beta-1}\prod_{j=1}^{r}\left(f_{j}t+g_{j}\right)^{\lambda_{j}}\Theta\left(\kappa_{l};-\frac{p}{t-a}-\frac{q}{b-t}\right)dt=\prod_{j=1}^{k}\left(af_{j}+g_{j}\right)^{\lambda_{j}}B\left(\alpha,\beta\right)\nonumber \\
\times\mathbf{F}_{D,\left(\kappa_{l}\right)}^{\left(r\right)}\left(\alpha,-\lambda_{1},\cdots,-\lambda_{r};\alpha+\beta;-\frac{\left(b-a\right)f_{1}}{af_{1}+g_{1}},\cdots,-\frac{\left(b-a\right)f_{r}}{af_{r}+g_{r}};\frac{p}{b-a},\frac{q}{b-a}\right)
\end{align}
where $a,b\in\mathbb{R}\ \left(a<b\right)$, $f_{i},g_{i},\lambda_{i}\in\mathbb{C},\ i=1,\cdots,k$,$\Re\left(b\right)>0,\Re\left(d\right)>0$,
$\Re\left(\alpha\right)>0,\ \Re\left(\beta\right)>0$ and
\[
\max\left\{ \left|\frac{\left(b-a\right)f_{1}}{af_{1}+g_{1}}\right|,\cdots,\left|\frac{\left(b-a\right)f_{r}}{af_{r}+g_{r}}\right|\right\} <1.
\]
\end{Theorem}
\begin{Proof}
Our proof mainly depends on the following expansion \cite{ Fractional Integration of the H-Function of Several Variables}
\begin{equation}
\left(f_{j}t+g_{j}\right)^{\lambda_{j}}=\left(af_{j}+g_{j}\right)^{\lambda_{j}}\sum_{m_{j}=0}^{\infty}\frac{\left(-\lambda_{j}\right)_{m_{j}}}{m_{j}!}\left(-\frac{\left(t-a\right)f_{j}}{af_{j}+g_{j}}\right)^{m_{j}}\ \left(\left|\left(t-a\right)f_{j}\right|<\left|af_{j}+g_{j}\right|;\ t\in\left[a,\ b\right]\right).
\end{equation}
Opening up $\prod_{j=1}^{r}\left(f_{j}t+g_{j}\right)^{\lambda_{j}}$
by using (72) and then integrating out $t$ we have
\begin{align}
 & \int_{a}^{b}\left(t-a\right)^{\alpha-1}\left(b-t\right)^{\beta-1}\prod_{j=1}^{r}\left(f_{j}t+g_{j}\right)^{\lambda_{j}}\Theta\left(\kappa_{l};-\frac{p}{t-a}-\frac{q}{b-t}\right)dt\nonumber \\
 & =\prod_{j=1}^{r}\left(af_{j}+g_{j}\right)^{\lambda_{j}}\sum_{m_{1},\cdots m_{r}=0}^{\infty}\frac{\left(-\lambda_{1}\right)_{m_{1}}\cdots\left(-\lambda_{r}\right)_{m_{r}}}{m_{1}!\cdots m_{r}!}\left(-\frac{f_{1}}{af_{1}+g_{1}}\right)^{m_{1}}\cdots\left(-\frac{f_{r}}{af_{r}+g_{r}}\right)^{m_{r}}\nonumber \\
 & \ \ \times\int_{a}^{b}\left(t-a\right)^{\alpha+m_{1}+\cdots+m_{r}-1}\left(b-t\right)^{\beta-1}\Theta\left(\kappa_{l};-\frac{p}{t-a}-\frac{q}{b-t}\right)dt.
\end{align}
The last integral in (73) can be evaluate by applying substitution $t=\left(b-a\right)u+a$, namely,
\begin{align}
 & \int_{a}^{b}\left(t-a\right)^{\alpha+m_{1}+\cdots+m_{r}-1}\left(b-t\right)^{\beta-1}\Theta\left(\kappa_{l};-\frac{p}{t-a}-\frac{q}{b-t}\right)dt\nonumber \\
 & =\left(b-a\right)^{\alpha+\beta+m_{1}+\cdots+m_{r}-1}\int_{0}^{1}u^{\alpha+m_{1}+\cdots+m_{r}-1}\left(1-u\right)^{\beta-1}\Theta\left(\kappa_{l};-\frac{p}{\left(b-a\right)u}-\frac{q}{\left(b-a\right)\left(1-u\right)}\right)du\nonumber \\
 & =\left(b-a\right)^{\alpha+\beta+m_{1}+\cdots+m_{r}-1}\mathcal{B}_{\frac{p}{b-a},\frac{q}{b-a}}^{\left(\kappa_{l}\right)}\left(\alpha+m_{1}+\cdots+m_{r},\beta\right).
\end{align}
Thus,
\begin{align}\nonumber
& \int_{a}^{b}\left(t-a\right)^{\alpha-1}\left(b-t\right)^{\beta-1}\prod_{j=1}^{r}\left(f_{j}t+g_{j}\right)^{\lambda_{j}}\Theta\left(\kappa_{l};-\frac{p}{t-a}-\frac{q}{b-t}\right)dt\\\nonumber
&
=\left(b-a\right)^{\alpha+\beta-1}\prod_{j=1}^{r}\left(af_{j}+g_{j}\right)^{\lambda_{j}}\nonumber
\sum_{m_{1},\cdots m_{r}=0}^{\infty}\left(-\lambda_{1}\right)_{m_{1}}\cdots\left(-\lambda_{r}\right)_{m_{r}}\mathcal{B}_{\frac{p}{b-a},\frac{q}{b-a}}^{\left(\kappa_{l}\right)}\left(\alpha+m_{1}+\cdots+m_{r},\beta\right)\\\nonumber
&
\times\frac{1}{m_{1}!}\left(-\frac{\left(b-a\right)f_{1}}{af_{1}+g_{1}}\right)^{m_{1}}\cdots\frac{1}{m_{r}!}\left(-\frac{\left(b-a\right)f_{r}}{af_{r}+g_{r}}\right)^{m_{r}}\\\nonumber
&
=\left(b-a\right)^{\alpha+\beta-1}\prod_{j=1}^{r}\left(af_{j}+g_{j}\right)^{\lambda_{j}}B\left(\alpha,\beta\right)\\\nonumber
&
\times\mathbf{F}_{D,\left(\kappa_{l}\right)}^{\left(r\right)}\left(\alpha,-\lambda_{1},\cdots,-\lambda_{r};\alpha+\beta;-\frac{\left(b-a\right)f_{1}}{af_{1}+g_{1}},\cdots,-\frac{\left(b-a\right)f_{r}}{af_{r}+g_{r}};\frac{p}{b-a},\frac{q}{b-a}\right).\nonumber
\end{align}
The second equality follows by the definition of extended Lauricella hypergeometric function $\mathbf{F}_{D,\left(\kappa_{l}\right)}^{\left(r\right)}$.
\end{Proof}
\begin{Theorem}
\begin{align}
 & \int_{0}^{\infty}\cdots\int_{0}^{\infty}\prod_{j=1}^{r}t_{j}^{\beta_{j}-1}e^{-t_{j}}{}_{1}F_{1}^{\left(\kappa_{l}\right)}\left[\begin{array}{c}
\alpha\\
\gamma
\end{array};x_{1}t_{1}+\cdots+x_{r}t_{r};b,d\right]dt_{1}\cdots dt_{r}\nonumber \\
 & \ \ =\Gamma\left(\beta_{1}\right)\cdots\Gamma\left(\beta_{r}\right)\mathbf{F}_{D,\left(\kappa_{l}\right)}^{\left(r\right)}\left(\alpha,\beta_{1},\cdots,\beta_{r};\gamma;x_{1},\cdots,x_{r};b,d\right).
\end{align}
\[
\left(\Re\left(\beta_1\right),\cdots,\Re\left(\beta_r\right)>0\right)
\]
\end{Theorem}
\begin{Proof}
The demonstration of this theorem depends on the application of the following series identity
\begin{equation}
\sum_{m_{1},\cdots,m_{r}=0}^{\infty}f\left(m_{1}+\cdots+m_{r}\right)\frac{x_{1}^{m_{1}}\cdots x_{r}^{m_{r}}}{m_{1}!\cdots m_{r}!}=\sum_{m=0}^{\infty}f\left(m\right)\frac{\left(x_{1}+\cdots+x_{r}\right)^{m}}{m!}.
\end{equation}
Specifically, we have
\begin{alignat*}{1}
 & \int_{0}^{\infty}\cdots\int_{0}^{\infty}\prod_{j=1}^{r}t_{j}^{\beta_{j}-1}e^{-t_{j}}{}_{1}F_{1}^{\left(\kappa_{l}\right)}\left(\alpha;\gamma;x_{1}t_{1}++x_{r}t_{r};b,d\right)dt_{1}\cdots dt_{r}\\
 & \ \ =\sum_{m_{1},\cdots,m_{r}=0}^{\infty}\frac{\mathcal{B}_{b,d}^{\left(\kappa_{l}\right)}\left(\alpha+m_{1}+\cdots+m_{r},\gamma-\alpha\right)}{B\left(\alpha,\gamma-\alpha\right)}\frac{x_{1}^{m_{1}}\cdots x_{r}^{m_{r}}}{m_{1}!\cdots m_{r}!}\prod_{j=1}^{r}\int_{0}^{\infty}t_{j}^{\beta_{j}+m_{j}-1}e^{-t_{j}}dt_{j}\\
 & \ \ =\sum_{m_{1},\cdots,m_{r}=0}^{\infty}\Gamma\left(\beta_{1}+m_{1}\right)\cdots\Gamma\left(\beta_{r}+m_{r}\right)\frac{\mathcal{B}_{b,d}^{\left(\kappa_{l}\right)}\left(\alpha+m_{1}+\cdots+m_{r},\gamma-\alpha\right)}{B\left(\alpha,\gamma-\alpha\right)}\frac{x_{1}^{m_{1}}\cdots x_{r}^{m_{r}}}{m_{1}!\cdots m_{r}!}\\
 & \ \ =\Gamma\left(\beta_{1}\right)\cdots\Gamma\left(\beta_{r}\right)\sum_{m_{1},\cdots,m_{r}=0}^{\infty}\frac{\mathcal{B}_{b,d}^{\left(\kappa_{l}\right)}\left(\alpha+m_{1}+\cdots+m_{r},\gamma-\alpha\right)}{B\left(\alpha,\gamma-\alpha\right)}\frac{x_{1}^{m_{1}}\cdots x_{r}^{m_{r}}}{m_{1}!\cdots m_{r}!}.
\end{alignat*}
\end{Proof}

In accordance with the previous method, we can obtain the following
generalization for classical Lauricella's hypergeometric function:
$\mathbf{F}_{A}^{\left(r\right)}\left(\alpha,\beta_{1},\cdots,\beta_{r};\gamma_{1},\cdots\gamma_{r};x_{1},\cdots,x_{r}\right)$.
\begin{Definition}
Lauricella hypergeometric function $\mathbf{F}_{A}^{\left(r\right)}\left(\alpha,\beta_{1},\cdots,\beta_{r};\gamma_{1},\cdots\gamma_{r};x_{1},\cdots,x_{r}\right)$ is defined by
\begin{multline}
\mathbf{F}_{A,\left(\kappa_{l}\right)}^{\left(r\right)}\left(\alpha,\beta_{1},\cdots,\beta_{r};\gamma_{1},\cdots\gamma_{r};x_{1},\cdots,x_{r};b,d\right)\\
=\sum_{m_{1},\cdots m_{r}=0}^{\infty}\left(\alpha\right)_{m_{1}+\cdots+m_{r}}\prod_{j=1}^{r}\frac{\mathcal{B}_{b,d}^{\left(\kappa_{l}\right)}\left(\beta_{j}+m_{j},\gamma_{j}-\beta_{j}\right)}{B\left(\beta_{j},\gamma_{j}-\beta_{j}\right)}\frac{x_{1}^{m_{1}}\cdots x_{r}^{m_{r}}}{m_{1}!\cdots m_{r}!}.
\end{multline}
\[
\left(\left|x_{1}\right|+\cdots+\left|x_{r}\right|<1;\ \min\left\{ \Re\left(b\right),\Re\left(d\right)\right\} \geq0\right)
\]
\end{Definition}

When $b=d=0$, function (77) reduces to usual Lauricella hypergeometric function \cite{ The H-Function Theory and Application}:
\begin{equation}
\mathbf{F}_{A}^{\left(r\right)}\left(\alpha,\beta_{1},\cdots,\beta_{r};\gamma_{1},\cdots\gamma_{r};x_{1},\cdots,x_{r}\right)
=\sum_{m_{1},\cdots m_{r}=0}^{\infty}\left(\alpha\right)_{m_{1}+\cdots+m_{r}}\frac{(b_1)_{m_1}\cdots(b_r)_{m_r}}{(c_1)_{m_1}\cdots(c_r)_{m_r}}\frac{x_{1}^{m_{1}}\cdots x_{r}^{m_{r}}}{m_{1}!\cdots m_{r}!}.\nonumber
\end{equation}
Comparing the coefficients of $\mathbf{F}_{A,\left(\kappa_{l}\right)}^{\left(r\right)}$ with those of $\mathbf{F}_{A}^{\left(r\right)}$, associated with \textbf{Definition 2.1}, it is clear that the construction of the coefficients of function (77) corresponds with what we generalize hypergeometric function.

For this new function, we have the following multidimensional integral representation.
\begin{Theorem} For Lauricella hypergeometric function $\mathbf{F}_{A}^{\left(r\right)}\left(\alpha,\beta_{1},\cdots,\beta_{r};\gamma_{1},\cdots\gamma_{r};x_{1},\cdots,x_{r}\right)$, we have
\begin{multline}
\mathbf{F}_{A,\left(\kappa_{l}\right)}^{\left(r\right)}\left(\alpha,\beta_{1},\cdots,\beta_{r};\gamma_{1},\cdots\gamma_{r};x_{1},\cdots,x_{r};b,d\right) =\prod_{j=1}^{r}B\left(\beta_{j},\gamma_{j}-\beta_{j}\right)\\
\times\int_{0}^{1}\cdots\int_{0}^{1}\prod_{j=1}^{r}u_{j}^{\beta_{j}-1}\left(1-u_{j}\right)^{\gamma_{j}-\beta_{j}-1}
\Theta\left(\kappa_{l};-\frac{b}{u_{j}}-\frac{d}{1-u_{j}}\right)\left(1-x_{1}u_{1}-\cdots-x_{r}u_{r}\right)^{-\alpha}du_{1}\cdots du_{r}.
\end{multline}
\[
\left(\Re\left(\beta_{j}\right)>0,\ \Re\left(\gamma_{j}-\beta_{j}\right)>0,\ j=1,\cdots,r;\ \min\left\{ \Re\left(b\right),\Re\left(d\right)\right\} \geq0\right)
\]
\end{Theorem}
\begin{Proof}This formula can be easily established by expanding the factor
$\left(1-x_{1}u_{1}-\cdots-x_{r}u_{r}\right)^{-\alpha}$ as
\[
\left(1-x_{1}u_{1}-\cdots-x_{r}u_{r}\right)^{-\alpha}=\sum_{m_{1},\cdots,m_{r}=0}^{\infty}\left(\alpha\right)_{m_{1}+\cdots+m_{r}}\frac{\left(x_{1}u_{1}\right)^{m_{1}}\cdots\left(x_{r}u_{r}\right)^{m_{r}}}{m_{1}!\cdots m_{r}!}
\]
and then integrating out $u_{j},\ j=1,\cdots,r$ with the help of our extended beta integral (11).
\end{Proof}

The following theorem show that extended Lauricella hypergeometric function $\mathbf{F}_{A,\left(\kappa_{l}\right)}^{\left(r\right)}$ can
be expressed as a single integral whose integrand is a product of several extended Kummer hypergeometric functions ${}_{1}F_{1}^{(\kappa_{l})}\left(\beta;\gamma;z;b,d\right)$.
\begin{Theorem}
\begin{align}
 & \mathbf{F}_{A,\left(\kappa_{l}\right)}^{\left(r\right)}\left(\alpha,\beta_{1},\cdots,\beta_{r};\gamma_{1},\cdots\gamma_{r};x_{1},\cdots,x_{r};b,d\right)\nonumber \\
 & \ \ =\frac{1}{\Gamma\left(\alpha\right)}\int_{0}^{1}e^{-t}t^{\alpha-1}{}_{1}F_{1}^{\left(\kappa_{l}\right)}\left[\begin{array}{c}
\beta_{1}\\
\gamma_{1}
\end{array};x_{1}t;b,d\right]\cdots{}_{1}F_{1}^{\left(\kappa_{l}\right)}\left[\begin{array}{c}
\beta_{r}\\
\gamma_{r}
\end{array};x_{r}t;b,d\right]dt.
\end{align}
\end{Theorem}
\begin{Proof}Taking the series form of $_{1}F_{1}^{\left(\kappa_{l}\right)}\left(\beta_{j};\gamma_{j};x_{j}t;b,d\right),\ j=1,\cdots,r$
and then integrating out $t$.
\end{Proof}

By applying the relation $\left(\alpha\right)_{m_{1}+\cdots+m_{r-1}+m_{r}}=\left(\alpha\right)_{m_{1}+\cdots+m_{r-1}}\left(\alpha+m_{1}+\cdots+m_{r-1}\right)_{m_{r}}$,
we can write
\begin{multline}
\mathbf{F}_{A,\left(\kappa_{l}\right)}^{\left(r\right)}\left(\alpha,\beta_{1},\cdots,\beta_{r};\gamma_{1},\cdots\gamma_{r};x_{1},\cdots,x_{r};b,d\right)\\ =\sum_{m_{1},\cdots m_{r-1}=0}^{\infty}\left(\alpha\right)_{m_{1}+\cdots+m_{r-1}}\prod_{j=1}^{r-1}\frac{\mathcal{B}_{b,d}^{\left(\kappa_{l}\right)}\left(\beta_{j}+m_{j},\gamma_{j}-\beta_{j}\right)}{B\left(\beta_{j},\gamma_{j}-\beta_{j}\right)}\frac{x_{1}^{m_{1}}\cdots x_{r-1}^{m_{r-1}}}{m_{1}!\cdots m_{r-1}!}\\
\times{}_{2}F_{1}^{\left(\kappa_{l}\right)}\left(\begin{array}{c}
\alpha+m_{1}+\cdots+m_{r-1},\beta_{r}\\
\gamma_{r}
\end{array};x_{r};b,d\right).
\end{multline}
\begin{Theorem}
The Mellin-Barnes type integral representation of $\mathbf{F}_{A,\left(\kappa_{l}\right)}^{\left(r\right)}$ is
\begin{align}
 & \mathbf{F}_{A,\left(\kappa_{l}\right)}^{\left(r\right)}\left(\alpha,\beta_{1},\cdots,\beta_{r};\gamma_{1},\cdots\gamma_{r};x_{1},\cdots,x_{r};b,d\right)\nonumber \\
 & \ \ =\frac{1}{\left(2\pi i\right)^{n}}\int_{L_{1}}\cdots\int_{L_{r}}\frac{\Gamma\left(\alpha-s_{1}-\cdots-s_{r}\right)}{\Gamma\left(\alpha\right)}\prod_{j=1}^{r}\frac{\mathcal{B}_{b,d}^{\left(\kappa_{l}\right)}\left(\beta_{j}-s_{j},\gamma_{j}-\beta_{j}\right)\Gamma\left(s_{j}\right)}{B\left(\beta_{j},\gamma_{j}-\beta_{j}\right)}\nonumber \\
 & \ \ \ \ \times\left(-x_{1}\right)^{-s_{1}}\cdots\left(-x_{r}\right)^{-s_{r}}ds_{1}\cdots ds_{r}.
\end{align}
where $L_{1},\cdots,L_{r}$ are Barnes paths of integration.
\end{Theorem}
\begin{Proof}
The basic ideas are the same when proving this theorem and \textbf{Theorem 3.4}. We take multidimensional Mellin transform of function $\mathbf{F}_{A,\left(\kappa_{l}\right)}^{\left(r\right)}\left(\alpha,\beta_{1},\cdots,\beta_{r};\gamma_{1},\cdots\gamma_{r};-x_{1},\cdots,-x_{r};b,d\right)$, get the result on the basis of applying the \emph{method of Brackets} and take the inverse Mellin transform of the function we get.

To obtain the multidimensional Mellin transform of $\mathbf{F}_{A,\left(\kappa_{l}\right)}^{\left(r\right)}\left(\alpha,\beta_{1},\cdots,\beta_{r};\gamma_{1},\cdots\gamma_{r};-x_{1},\cdots,-x_{r};b,d\right)$,
we multiply both sides of (77) by $x_{1}^{s_{1}-1}\cdots x_{r}^{s_{r}-1}$ and integrate with respect to $x_{1},\cdots,x_{r}$ over $\left[0,\infty\right)\times\cdots\times\left[0,\infty\right)$. One has
\begin{align*}
 & \int_{0}^{\infty}\cdots\int_{0}^{\infty}x_{1}^{s_{1}-1}\cdots x_{r}^{s_{r}-1}\mathbf{F}_{A,\left(\kappa_{l}\right)}^{\left(r\right)}\left(\alpha,\beta_{1},\cdots,\beta_{r};\gamma_{1},\cdots\gamma_{r};-x_{1},\cdots,-x_{r};b,d\right)dx_{1}\cdots dx_{r}\\
 & =\sum_{m_{1},\cdots m_{r}=0}^{\infty}\frac{\Gamma\left(\alpha+m_{1}+\cdots+m_{r}\right)}{\Gamma\left(\alpha\right)}\prod_{j=1}^{r}\frac{\mathcal{B}_{b,d}^{\left(\kappa_{l}\right)}\left(\beta_{j}+m_{j},\gamma_{j}-\beta_{j}\right)}{B\left(\beta_{j},\gamma_{j}-\beta_{j}\right)}\frac{\left(-1\right)^{m_{1}}\cdots\left(-1\right)^{m_{r}}}{m_{1}!\cdots m_{r}!}\\
 & \ \ \times\int_{0}^{\infty}\cdots\int_{0}^{\infty}x_{1}^{m_{1}+s_{1}-1}\cdots x_{r}^{m_{r}+s_{1}-1}dx_{1}\cdots dx_{r}\\
 & =\sum_{m_{1},\cdots m_{r}=0}^{\infty}f\left(m_{1},\cdots,m_{r}\right)\phi_{m_{1}}\cdots\phi_{m_{r}}\left\langle s_{1}+m_{1}\right\rangle \cdots\left\langle s_{r}+m_{r}\right\rangle.
\end{align*}
The result here is just more complicated in form than that of equation (56). And with the basic ideas when we prove \textbf{Theorem 3.4} and with the applying of \textbf{Rule III}\cite{Ramanujan Master Theorem}, we can get
\begin{align*}
\int_{0}^{\infty}\cdots\int_{0}^{\infty}x_{1}^{s_{1}-1}\cdots x_{r}^{s_{r}-1}\mathbf{F}_{A,\left(\kappa_{l}\right)}^{\left(r\right)}\left(\alpha,\beta_{1},\cdots,\beta_{r};\gamma_{1},\cdots\gamma_{r};-x_{1},\cdots,-x_{r};b,d\right)dx_{1}\cdots dx_{r}\\
=\frac{\Gamma\left(\alpha-s_{1}-\cdots-s_{r}\right)}{\Gamma\left(\alpha\right)}\prod_{j=1}^{r}\frac{\mathcal{B}_{b,d}^{\left(\kappa_{l}\right)}\left(\beta_{j}-s_{j},\gamma_{j}-\beta_{j}\right)\Gamma\left(s_{j}\right)}{B\left(\beta_{j},\gamma_{j}-\beta_{j}\right)}.
\end{align*}
Taking multidimensional inverse Mellin transform we have
\begin{align*}
 & \mathbf{F}_{A,\left(\kappa_{l}\right)}^{\left(r\right)}\left(\alpha,\beta_{1},\cdots,\beta_{r};\gamma_{1},\cdots\gamma_{r};-x_{1},\cdots,-x_{r};b,d\right)\\
 & =\frac{1}{\left(2\pi i\right)^{n}}\int_{L_{1}}\cdots\int_{L_{r}}\frac{\Gamma\left(\alpha-s_{1}-\cdots-s_{r}\right)}{\Gamma\left(\alpha\right)}\prod_{j=1}^{r}\frac{\mathcal{B}_{b,d}^{\left(\kappa_{l}\right)}\left(\beta_{j}-s_{j},\gamma_{j}-\beta_{j}\right)\Gamma\left(s_{j}\right)}{B\left(\beta_{j},\gamma_{j}-\beta_{j}\right)}\\
 & \ \ x_{1}^{-s_{1}}\cdots x_{r}^{-s_{r}}ds_{1}\cdots ds_{r}.
\end{align*}
Theorem 3.13 follows from replacing $-x_i, i=1,\cdots,r$ with $x_i$.
\end{Proof}

The way of proving the following two theorems are exactly the same.
\begin{Theorem}
The Mellin-Barnes type integral representation of  $\mathbf{F}_{D,\left(\kappa_{l}\right)}^{\left(r\right)}$ is
\begin{multline}
\mathbf{F}_{D,\left(\kappa_{l}\right)}^{\left(r\right)}\left(\alpha,\beta_{1},\cdots,\beta_{r};\gamma;x_{1},\cdots,x_{r};b,d\right) =\frac{1}{\left(2\pi i\right)^{n}}\\
\times\int_{L_{1}}\cdots\int_{L_{r}}\prod_{j=1}^{r}\frac{\Gamma\left(\beta_{j}-s_{j}\right)\Gamma\left(s_{j}\right)}{\Gamma\left(\beta_{j}\right)}\frac{\mathcal{B}_{b,d}^{\left(\kappa_{l}\right)}\left(\alpha-s_{1}-\cdots-s_{r},\gamma-\alpha\right)}{B\left(\alpha,\gamma-\alpha\right)}\left(-x_{1}\right)^{-s_{1}}\cdots\left(-x_{r}\right)^{-s_{2}}ds_{1}\cdots ds_{r}
\end{multline}
where $L_{1},\cdots,L_{r}$ are Barnes paths of integration.
\end{Theorem}
\begin{Theorem}
For $p=q+1$, we have
\begin{align}
 & \Gamma\left(\alpha_{1}\right)\prod_{j=1}^{q}\frac{\Gamma\left(\alpha_{j+1}\right)\Gamma\left(\beta_{j}-\alpha_{j+1}\right)}{\Gamma\left(\beta_{j}\right)}{}_{q+1}F_{q}^{\left(\kappa_{l}\right)}\left(\begin{array}{ccc}
\alpha_{1} & \cdots & \alpha_{q+1}\\
\beta_{1} & \cdots & \beta_{q}
\end{array};-\left(x_{1}+\cdots+x_{r}\right);b,d\right)\nonumber \\
 & \ \ =\frac{1}{\left(2\pi i\right)^{n}}\int_{L_{1}}\cdots\int_{L_{r}}\Gamma\left(\alpha_{1}-s_{1}-\cdots-s_{r}\right)\prod_{j=1}^{q}\mathcal{B}_{b,d}^{\left(\kappa_{l}\right)}\left(\alpha_{j+1}-s_{1}-\cdots-s_{r},\beta_{j}-\alpha_{j+1}\right)\nonumber \\
 & \ \ \times\Gamma\left(s_{1}\right)\cdots\Gamma\left(s_{r}\right)x_{1}^{-s_{1}}\cdots x_{r}^{-s_{r}}ds_{1}\cdots ds_{r}
\end{align}
For $p=q$
\begin{align}
 & \prod_{j=1}^{q}\frac{\Gamma\left(\alpha_{j}\right)\Gamma\left(\beta_{j}-\alpha_{j}\right)}{\Gamma\left(\beta_{j}\right)}{}_{q}F_{q}^{\left(\kappa_{l}\right)}\left(\begin{array}{ccc}
\alpha_{1} & \cdots & \alpha_{q}\\
\beta_{1} & \cdots & \beta_{q}
\end{array};-\left(x_{1}+\cdots+x_{r}\right);b,d\right)\nonumber \\
 & \ \ =\frac{1}{\left(2\pi i\right)^{n}}\int_{L_{r}}\cdots\int_{L_{r}}\prod_{j=1}^{q}\mathcal{B}_{b,d}^{\left(\kappa_{l}\right)}\left(\alpha_{j}-s_{1}-\cdots-s_{r},\beta_{j}-\alpha_{j}\right)\nonumber \\
 & \ \ \times\Gamma\left(s_{1}\right)\cdots\Gamma\left(s_{r}\right)x_{1}^{-s_{1}}\cdots x_{r}^{-s_{r}}ds_{1}\cdots ds_{r}
\end{align}
For $p<q$
\begin{align}
 & \frac{1}{\Gamma\left(\beta_{1}\right)\cdots\Gamma\left(\beta_{r}\right)}\prod_{j=1}^{p}\frac{\Gamma\left(\alpha_{j}\right)\Gamma\left(\beta_{r+j}-\alpha_{j}\right)}{\Gamma\left(\beta_{r+j}\right)}{}_{p}F_{q}^{\left(\kappa_{l}\right)}\left(\begin{array}{ccc}
\alpha_{1} & \cdots & \alpha_{p}\\
\beta_{1} & \cdots & \beta_{q}
\end{array};-\left(x_{1}+\cdots+x_{r}\right);b,d\right)\nonumber \\
 & \ \ =\frac{1}{\left(2\pi i\right)^{n}}\int_{L_{1}}\cdots\int_{L_{r}}\prod_{j=1}^{r}\frac{1}{\Gamma\left(\beta_{j}-s_{1}-\cdots-s_{r}\right)}\prod_{j=1}^{p}\mathcal{B}_{b,d}^{\left(\kappa_{l}\right)}\left(\alpha_{j}-s_{1}-\cdots-s_{r},\beta_{r+j}-\alpha_{j}\right)\nonumber \\
 & \ \ \times\Gamma\left(s_{1}\right)\cdots\Gamma\left(s_{r}\right)x_{1}^{-s_{1}}\cdots x_{r}^{-s_{r}}ds_{1}\cdots ds_{r}
\end{align}
where $L_{1},\cdots,L_{r}$ are Barnes paths of integration.
\end{Theorem}
\begin{Remark}Note that when parameters $b=d=0$, integral (83), (84) and (85)
reduce to
\begin{align}
 & \frac{\prod_{j=1}^{p}\Gamma\left(\alpha_{j}\right)}{\prod_{j=1}^{q}\Gamma\left(\beta_{j}\right)}{}_{p}F_{q}\left[\begin{array}{ccc}
\alpha_{1} & \cdots & \alpha_{p}\\
\beta_{1} & \cdots & \beta_{q}
\end{array};-\left(x_{1}+\cdots+x_{r}\right)\right]\nonumber \\
 & \ \ =\frac{1}{\left(2\pi i\right)^{n}}\int_{L_{1}}\cdots\int_{L_{r}}\frac{\prod_{j=1}^{p}\Gamma\left(\alpha_{j}+s_{1}+\cdots+s_{r}\right)}{\prod_{j=1}^{q}\Gamma\left(\beta_{j}+s_{1}+\cdots+s_{r}\right)}\Gamma\left(-s_{1}\right)\cdots\Gamma\left(-s_{r}\right)x_{1}^{s_{1}}\cdots x_{r}^{s_{r}}ds_{1}\cdots ds_{r}
\end{align}
where the contours are Barnes type with indentations, if necessary,
such that the poles of $\Gamma\left(\alpha_{j}+s_{1}+\cdots+s_{r}\right),\ j=1,\cdots,p$ are separated from those of $\Gamma\left(-s_{j}\right),\ j=1,\cdots,r$.
And this multiple integral representation for the generalized hypergeometric
series is given by Saigo and Saxena \cite[p. 212]{The H-Function Theory and Application}.
\end{Remark}

\section{{\normalsize Applications in Hilbert-Hardy type Inequalities}}

We have introduced many different generalizations of usual hypergeometric
functions, such as, extended generalized hypergeometric function,
Appell's and Lauricella's hypergeometric functions. A lot of important
results are established in the previous sections. Now we hope to find
some applications in other branches of mathematics. In this section,
we will use extended Gauss hypergeometric function to establish several
important inequalities.

Let $f\left(x\right),g\left(y\right)\geq0$, $f\left(x\right)\in L^{p}\left(0,+\infty\right),g\left(x\right)\in L^{q}\left(0,+\infty\right)$,
$\frac{1}{p}+\frac{1}{q}=1$ and $p>1$. Then we have the following
two equivalent inequalities as
\begin{equation}
\int_{0}^{\infty}\int_{0}^{\infty}\frac{f\left(x\right)g\left(y\right)}{x+y}dxdy\leq\frac{\pi}{\sin\left(\frac{\pi}{p}\right)}\left\{ \int_{0}^{\infty}f^{p}\left(x\right)dx\right\} ^{\frac{1}{p}}\left\{ \int_{0}^{\infty}g^{q}\left(y\right)dy\right\} ^{\frac{1}{q}}
\end{equation}
\begin{equation}
\int_{0}^{\infty}\left(\int_{0}^{\infty}\frac{f\left(x\right)}{x+y}dx\right)^{p}dy\leq\left[\frac{\pi}{\sin\left(\frac{\pi}{p}\right)}\right]^{p}\int_{0}^{\infty}f^{p}\left(x\right)dx
\end{equation}
where the constant factor $\frac{\pi}{\sin\left(\frac{\pi}{p}\right)}$
and $\left[\frac{\pi}{\sin\left(\frac{\pi}{p}\right)}\right]^{p}$
are the best possible in (87) and (88) respectively. And the equality
in (87) and (88) holds iff $f\left(x\right)=0$ or $g\left(x\right)=0$.
Inequality (87) is called Hardy-Hilbert integral inequality (see \cite{Some Extensions of Hilberts Integral Inequality}, \cite{Hilbert Inequality and Gaussian Hypergeometric Functions}).
During the last decade inequality (87) and (88) were generalized in many
different ways, some of the generalizations rely heavily on the integral
identities involving hypergeometric functions and beta functions (see \cite{Some Extensions of Hilberts Integral Inequality}, \cite{Hilbert Inequality and Gaussian Hypergeometric Functions}). In \cite{Hilbert Inequality and Gaussian Hypergeometric Functions}, the author obtain some new Hardy-Hilbert
type inequalities with fractional kernels, which involve the constants
expressed in terms of Gauss hypergeometric functions. So it is natural
to consider whether these new inequalities can be extended to the
more general case by using our extended Gauss hypergeometric functions. The answer is positive. In the whole process of proving, we can get some previously unimagined properties. Note that in what follows, we chose function: $\Theta\left(\kappa_{l};z\right)=\exp\left(z\right)$. This expression are easy to handle and the final form of our result will thus be more simple.

We need the following two integral identities.
\begin{Lemma} Suppose $a,b,c,\alpha,\gamma,\tilde p,\tilde q\in\mathbb{R}$ are
s.t. $\tilde p\geq0, \tilde q\geq0, \ a+c>b>0$ and $0<\alpha<2\gamma$. Then
\begin{equation}
\int_{0}^{\infty}\frac{x^{b-1}}{\left(1+\gamma x\right)^{c}\left(1+\alpha x\right)^{a}}\exp\left(-\gamma\tilde{q}x-\frac{\tilde{p}\gamma^{-1}}{x}\right)dx=e^{\tilde{p}+\tilde{q}}\gamma^{-b}B\left(b,c+a-b\right){}_{2}F_{1}\left[\begin{array}{c}
a,b\\
c+a
\end{array};\frac{\gamma-\alpha}{\gamma};\tilde{p},\tilde{q}\right]
\end{equation}
\begin{equation}
\int_{0}^{\infty}\frac{x^{b-1}}{\left(1+\alpha x\right)^{a}\left(1+\gamma x\right)^{c}}\exp\left(-\alpha\tilde{q}x-\frac{\tilde{p}\alpha^{-1}}{x}\right)dx=e^{\tilde{p}+\tilde{q}}\alpha^{-b}B\left(b,c+a-b\right){}_{2}F_{1}\left[\begin{array}{c}
c,b\\
c+a
\end{array};\frac{\alpha-\gamma}{\alpha};\tilde{p},\tilde{q}\right]
\end{equation}
\end{Lemma}
\begin{Proof} We start with the integral representation (13) with $\Theta\left(\kappa_{l};z\right)=\exp\left(z\right)$, i.e.,
\begin{equation}
\int_{0}^{1}t^{b-1}\left(1-t\right)^{c-b-1}\left(1-zt\right)^{-a}\exp\left(-\frac{\tilde{p}}{t}-\frac{\tilde{q}}{1-t}\right)dt=B\left(b,c-b\right){}_{2}F_{1}\left[\begin{array}{c}
a,b\\
c
\end{array};z;\tilde{p},\tilde{q}.\right]
\end{equation}
By using the substitutions $1-t=\frac{1}{1+u}$ in (91), we have
\begin{equation}
\int_{0}^{\infty}u^{b-1}\left(1+u\right)^{a-c}\left(1+u-zu\right)^{-a}\exp\left(-u\tilde{q}-\frac{\tilde{p}}{u}\right)du=e^{\tilde{p}+\tilde{q}}B\left(b,c-b\right){}_{2}F_{1}\left[\begin{array}{c}
a,b\\
c
\end{array};z;\tilde{p},\tilde{q}\right]
\end{equation}
Then we substitute $\gamma x$ for $u$ in Eq (92)
\begin{equation}
\int_{0}^{\infty}\frac{x^{b-1}}{\left(1+\gamma x\right)^{c-a}\left(1+\left(\gamma-z\gamma\right)x\right)^{a}}\exp\left(-\gamma\tilde{q}x-\frac{\tilde{p}\gamma^{-1}}{x}\right)du=e^{\tilde{p}+\tilde{q}}\gamma^{-b}B\left(b,c-b\right){}_{2}F_{1}\left[\begin{array}{c}
a,b\\
c
\end{array};z;\tilde{p},\tilde{q}\right]
\end{equation}
The first assertion of \textbf{lemma 2} follows from setting $\gamma\left(1-z\right)=\alpha$
and replace $c-a$ with $c$ in above equation.\\
To prove the second assertion we can make transformation $u=\alpha x$
in Eq. (92) and then set $\alpha\left(1-z\right)=\gamma$. Replacing
$c-a$ with $c$ in the result and then interchange $a$ and $c$
we obtain the relation (90).
\end{Proof}
\begin{Remark}
If we let $\tilde{q}=\tilde{p}=0$, \textbf{lemma 2} reduce to \cite[LEMMA
1.]{Hilbert Inequality and Gaussian Hypergeometric Functions}
\[
\int_{0}^{\infty}\frac{x^{b-1}}{\left(1+\alpha x\right)^{a}\left(1+\gamma x\right)^{c}}dx=\gamma^{-b}B\left(b,a+c-b\right)F\left[\begin{array}{c}
a,b\\
a+c
\end{array};\frac{\gamma-\alpha}{\gamma}\right].
\]
In fact, introducing the factor $\exp\left(-\frac{\tilde{p}}{t}-\frac{\tilde{q}}{1-t}\right)$
weakens the symmetry of integrand. So we need two integral identities
to establish our results.
\end{Remark}

It is time to finish the proof of \textbf{Theorem 2.5}.

\textbf{The Proof of Theorem 2.5:} With the help of the formulas (89) and (90) of \textbf{Lemma 2}
we are able to finish the proof of the Euler transformation mentioned
in Theorem 2.5. Let us set $\tilde{q}=\frac{\gamma\tilde{q}}{\alpha}$
and $\tilde{p}=\frac{\alpha\tilde{p}}{\gamma}$ in the second assertion
(90). Then (90) becomes
\begin{equation}
\int_{0}^{\infty}\frac{x^{b-1}}{\left(1+\alpha x\right)^{a}\left(1+\gamma x\right)^{c}}\exp\left(-\gamma\tilde{q}x-\frac{\tilde{p}\gamma^{-1}}{x}\right)dx=e^{\frac{\alpha\tilde{p}}{\gamma}+\frac{\gamma\tilde{q}}{\alpha}}\alpha^{-b}B\left(b,c+a-b\right){}_{2}F_{1}\left[\begin{array}{c}
c,b\\
c+a
\end{array};\frac{\alpha-\gamma}{\alpha};\frac{\alpha\tilde{p}}{\gamma},\frac{\gamma\tilde{q}}{\alpha}\right].
\end{equation}
By equalizing the identities (94) and (89) we obtain the relation:
\begin{equation}
e^{-\left(\frac{\frac{\gamma}{\alpha}-1}{\frac{\gamma}{\alpha}}\right)\tilde{p}-\left(1-\frac{\gamma}{\alpha}\right)\tilde{q}}\left(\frac{\gamma}{\alpha}\right)^{b}{}_{2}F_{1}\left[\begin{array}{c}
c,b\\
c+a
\end{array};1-\frac{\gamma}{\alpha};\frac{\alpha\tilde{p}}{\gamma},\frac{\gamma\tilde{q}}{\alpha}\right]={}_{2}F_{1}\left[\begin{array}{c}
a,b\\
c+a
\end{array};\frac{\frac{\gamma}{\alpha}-1}{\frac{\gamma}{\alpha}};\tilde{p},\tilde{q}\right].
\end{equation}
If we write $\frac{\gamma}{\alpha}=1-z$, then (93) becomes
\begin{equation}
e^{-\left(1-z\right)\tilde{p}-z\tilde{q}}\left(1-z\right)^{b}{}_{2}F_{1}\left[\begin{array}{c}
c-a,b\\
c
\end{array};z;\frac{\tilde{p}}{1-z},\left(1-z\right)\tilde{q}\right]={}_{2}F_{1}\left[\begin{array}{c}
a,b\\
c
\end{array};\frac{-z}{1-z};\tilde{p},\tilde{q}\right].
\end{equation}
By using the notation occured in \textbf{section 2}. we have:
\begin{equation}
e^{-\left(1-z\right)b-zd}\left(1-z\right)^{\alpha_{2}}{}_{2}F_{1}\left[\begin{array}{c}
\beta_{1}-\alpha_{1},\alpha_{2}\\
\beta_{1}
\end{array};z;\frac{b}{1-z},\left(1-z\right)d\right]={}_{2}F_{1}\left[\begin{array}{c}
\alpha_{1},\alpha_{2}\\
\beta_{1}
\end{array};\frac{-z}{1-z};b,d\right]
\end{equation}
Substituing equation (95) into (24) we get:
\begin{align}
_{2}F_{1}\left[\begin{array}{c}
\alpha_{1},\alpha_{2}\\
\beta_{1}
\end{array};z;b,d\right] & =\left(1-z\right)^{-\alpha_{1}}{}_{2}F_{1}\left[\begin{array}{c}
\alpha_{1},\beta_{1}-\alpha_{2}\\
\alpha_{2}
\end{array};\frac{-z}{1-z};d,b\right]\nonumber \\
 & =e^{-\left(1-z\right)b-zd}\left(1-z\right)^{\beta_{1}-\alpha_{2}-\alpha_{1}}{}_{2}F_{1}\left[\begin{array}{c}
\beta_{1}-\alpha_{1},\beta_{1}-\alpha_{2}\\
\beta_{1}
\end{array};z;\frac{b}{1-z},\left(1-z\right)d\right]
\end{align}
This is in fact a generalized Euler transformation.

Let us come back to the proof of inequalities. Before we state and prove our new inequality, we need to define some weighted functions. We define $F:\left(0,\infty\right)\mapsto\mathbb{R}$ by:
\begin{equation}
F\left(x;\tilde{p},\tilde{q}\right)=\left[\int_{0}^{\infty}\frac{y^{-q'A_{2}}}{\left(x+\alpha_{1}y\right)^{s_{1}}\left(x+\alpha_{2}y\right)^{s_{2}}}\exp\left(-\alpha_{1}\tilde{q}\frac{y}{x}-\frac{\tilde{p}}{\alpha_{1}}\frac{x}{y}\right)dy\right]^{\frac{1}{q'}},
\end{equation}
where $\tilde{p}\geq0,\ \tilde{q}\geq0$, $\mu>0,\ \alpha_{1},\alpha_{2}>0$,
$\frac{1}{2}<\frac{\alpha_{1}}{\alpha_{2}}<2$, $s_{1}+s_{2}>0$ and
$A_{2}\in\left(\frac{1-s_{1}-s_{2}}{q'},\frac{1}{q'}\right)$. Define $G:\left(0,\infty\right)\mapsto\mathbb{R}$ by:
\begin{equation}
G\left(y;\tilde{p},\tilde{q}\right)=\left[\int_{0}^{\infty}\frac{x^{-p'A_{1}}}{\left(x+\alpha_{1}y\right)^{s_{1}}\left(x+\alpha_{2}y\right)^{s_{2}}}\exp\left(-\frac{\tilde{q}}{\alpha_{2}}\frac{x}{y}-\frac{\tilde{p}}{\alpha_{2}}\frac{y}{x}\right)dx\right]^{\frac{1}{p'}},
\end{equation}
where $\tilde{p}\geq0,\ \tilde{q}\geq0$, $\mu>0,\alpha_{1},\alpha_{2}>0$,
$\frac{1}{2}<\frac{\alpha_{1}}{\alpha_{2}}<2$, $s_{1}+s_{2}>0$ and
$A_{1}\in\left(\frac{1-s_{1}-s_{2}}{p'},\frac{1}{p'}\right)$.
By virtue of \textbf{lemma 1} we can compute $F\left(x\right)$ and $G\left(y\right)$ and
express them as extended Gauss hypergeometric functions.
\begin{Lemma} Suppose $\mu>0,\alpha_{1},\alpha_{2}>0$, $\frac{1}{2}<\frac{\alpha_{1}}{\alpha_{2}}<2$,
$s_{1}+s_{2}>0$, Further, let $A_{1}$ and $A_{2}$ be real parameters
such that $A_{1}\in\left(\frac{1-s_{1}-s_{2}}{p'},\frac{1}{p'}\right)$
and $A_{2}\in\left(\frac{1-s_{1}-s_{2}}{q'},\frac{1}{q'}\right)$.
If the functions $F$ and $G$ are defined by (99) and (100) respectively,
then
\begin{equation}
F\left(x;\tilde{p},\tilde{q}\right)=e^{\frac{\tilde{p}+\tilde{q}}{q'}}k\left(F;\ \tilde{p},\tilde{q}\right)x^{\frac{1-s_{1}-s_{2}}{q'}-A_{2}},
\end{equation}
\begin{equation}
G\left(y;\tilde{p},\tilde{q}\right)=e^{\frac{\tilde{p}+\tilde{q}}{p'}}k\left(G;\ \tilde{p},\tilde{q}\right)y^{\frac{1-s_{1}-s_{2}}{p'}-A_{1}},
\end{equation}
where
\begin{equation}
k\left(F;\ \tilde{p},\tilde{q}\right)=\alpha_{1}^{A_{2}-\frac{1}{q'}}B^{\frac{1}{q'}}\left(1-q'A_{2},s_{1}+s_{2}+q'A_{2}-1\right)\left\{ _{2}F_{1}\left[\begin{array}{c}
s_{2},1-q'A_{2}\\
s_{1}+s_{2}
\end{array};\frac{\alpha_{1}-\alpha_{2}}{\alpha_{1}};\tilde{p},\tilde{q}\right]\right\} ^{\frac{1}{q'}},
\end{equation}
\begin{equation}
k\left(G;\ \tilde{p},\tilde{q}\right)=\alpha_{1}^{-\frac{s_{1}}{p'}}\alpha_{2}^{\frac{1-s_{2}}{p'}-A_{1}}B^{\frac{1}{p'}}\left(1-p'A_{1},s_{1}+s_{2}+p'A_{1}-1\right)\left\{ _{2}F_{1}\left[\begin{array}{c}
s_{1},1-p'A_{1}\\
s_{1}+s_{2}
\end{array};\frac{\alpha_{1}-\alpha_{2}}{\alpha_{1}};\tilde{p},\tilde{q}\right]\right\} ^{\frac{1}{p'}}.
\end{equation}
\end{Lemma}
\begin{Proof}This lemma is a direct consequence of \textbf{lemma 2}.
\end{Proof}

Now we are prepared to derive our main result.
\begin{Theorem}
Let $p$ and $q$ be the real parameters such that
$p>1,\ q>1,\ \frac{1}{p}+\frac{1}{q}\geq1,$and let $p'$ and $q'$
respectively be their conjugate exponents, that is, $\frac{1}{p}+\frac{1}{p'}=1$
and $\frac{1}{q}+\frac{1}{q'}=1$. Define $\lambda=\frac{1}{p'}+\frac{1}{q'}$
and note that $0<\lambda\leq1$, for all $p$ and $q$. Further, suppose
$\mu>0,\alpha_{1},\alpha_{2}>0,\frac{1}{2}<\frac{\alpha_{1}}{\alpha_{2}}<2,s_{1}+s_{2}>0$.
If $f$ and $g$ are non-negative measurable functions on $\left(0,\infty\right)$,
then the following inequalities:
\begin{multline}
\int_{0}^{\infty}\int_{0}^{\infty}\frac{\exp\left[-\left(\alpha_{1}\tilde{q}+\alpha_{2}\tilde{p}\right)\left(\frac{y}{x}+\frac{1}{\alpha_{1}\alpha_{2}}\frac{x}{y}\right)\right]}{\left(x+\alpha_{1}y\right)^{\lambda s_{1}}\left(x+\alpha_{2}y\right)^{\lambda s_{2}}}f\left(x\right)g\left(y\right)dxdy\\
 \leq
 e^{2\left(\tilde{p}+\tilde{q}\right)}K_{\tilde{p}.\tilde{q}}\left(\int_{0}^{\infty}x^{\frac{p}{q'}\left(1-s_{1}-s_{2}\right)+p\left(A_{1}-A_{2}\right)}f^{p}\left(x\right)dx\right)^{\frac{1}{p}}
\left(\int_{0}^{\infty}y^{\frac{q}{p'}\left(1-s_{1}-s_{2}\right)+q\left(A_{2}-A_{1}\right)}g^{q}\left(y\right)dy\right)^{\frac{1}{q}}
\end{multline}
and
\begin{align}
\left[\int_{0}^{\infty}y^{\frac{q'}{p'}\left(s_{1}+s_{2}-1\right)+q'\left(A_{1}-A_{2}\right)}\left(\int_{0}^{\infty}\frac{\exp\left[-\left(\alpha_{1}\tilde{q}+\alpha_{2}\tilde{p}\right)\left(\frac{y}{x}+\frac{1}{\alpha_{1}\alpha_{2}}\frac{x}{y}\right)\right]}{\left(x+\alpha_{1}y\right)^{\lambda s_{1}}\left(x+\alpha_{2}y\right)^{\lambda s_{2}}}f\left(x\right)dx\right)^{q'}dy\right]^{\frac{1}{q'}}\nonumber \\
\leq e^{2\left(\tilde{p}+\tilde{q}\right)}K_{\tilde{p}.\tilde{q}}\left(\int_{0}^{\infty}x^{\frac{p}{q'}\left(1-s_{1}-s_{2}\right)+p\left(A_{1}-A_{2}\right)}f^{p}\left(x\right)dx\right)^{\frac{1}{p}}
\end{align}
hold for any $A_{1}\in\left(\frac{1-s_{1}-s_{2}}{p'},\frac{1}{p'}\right)$
and $A_{2}\in\left(\frac{1-s_{1}-s_{2}}{q'},\frac{1}{q'}\right)$.
The constant
\[
K_{\tilde{p}.\tilde{q}}=k\left(F;\ q'\tilde{p},q'\tilde{q}\right)k\left(G;\ p'\tilde{p},p'\tilde{q}\right).
\]
Inequalities (105) and (106) are equivalent. In addition, the equality in
(105) and (106) holds iff at least one of the functions $f$ or $g$ is
equal to zero.
\end{Theorem}
\begin{Proof}The left-hand side of (105) can be written as
\begin{equation}
\int_{0}^{\infty}\int_{0}^{\infty}\frac{\exp\left[-\left(\alpha_{1}\tilde{q}+\alpha_{2}\tilde{p}\right)\left(\frac{y}{x}+\frac{1}{\alpha_{1}\alpha_{2}}\frac{x}{y}\right)\right]}{\left(x+\alpha_{1}y\right)^{\lambda s_{1}}\left(x+\alpha_{2}y\right)^{\lambda s_{2}}}f\left(x\right)g\left(y\right)dxdy=\int_{0}^{\infty}\int_{0}^{\infty}K_{1}^{\frac{1}{q'}}\left(x,y\right)K_{2}^{\frac{1}{p'}}\left(x,y\right)K_{3}^{1-\lambda}\left(x,y\right)dxdy
\end{equation}
where
\[
K_{1}\left(x,y\right)=\underset{I_{1}\left(x,y\right)}{\underbrace{\frac{F^{p-q'}\left(x;q'\tilde{p},q'\tilde{q}\right)f^{p}\left(x\right)}{\left(x+\alpha_{1}y\right)^{s_{1}}\left(x+\alpha_{2}y\right)^{s_{2}}}\frac{x^{pA_{1}}}{y^{q'A_{2}}}}}\exp\left(-\alpha_{1}q'\tilde{q}\frac{y}{x}-q'\frac{\tilde{p}\alpha_{1}^{-1}}{yx^{-1}}\right)
\]
\[
K_{2}\left(x,y\right)=\underset{I_{2}\left(x,y\right)}{\underbrace{\frac{G^{q-p'}\left(y;p'\tilde{p},p'\tilde{q}\right)g^{q}\left(y\right)}{\left(x+\alpha_{1}y\right)^{s_{1}}\left(x+\alpha_{2}y\right)^{s_{2}}}\frac{y^{qA_{2}}}{x^{p'A_{1}}}}}\exp\left(-\alpha_{2}^{-1}p'\tilde{q}\frac{x}{y}-p'\frac{\tilde{p}\alpha_{2}}{xy^{-1}}\right)
\]
\[
K_{3}\left(x,y\right)=x^{pA_{1}}F^{p}\left(x;q'\tilde{p},q'\tilde{q}\right)y^{qA_{2}}G^{q}\left(y;p'\tilde{p},p'\tilde{q}\right)f^{p}\left(x\right)g^{q}\left(y\right).
\]
(It is easy but time-consuming to verify Eq (107).) and the functions
$F$ and $G$ are defined by (99), (100). Now, since the exponents satisfy
identity: $\frac{1}{p'}+\frac{1}{q'}+\left(1-\lambda\right)=1$. We
can use Holder's inequality, namely,
\begin{align*}
 & \int_{0}^{\infty}\int_{0}^{\infty}K_{1}^{\frac{1}{q'}}\left(x,y\right)K_{2}^{\frac{1}{p'}}\left(x,y\right)K_{3}^{1-\lambda}\left(x,y\right)dxdy\\
 & \leq\left(\int_{0}^{\infty}\int_{0}^{\infty}K_{1}\left(x,y\right)dxdy\right)^{\frac{1}{q'}}\left(\int_{0}^{\infty}\int_{0}^{\infty}K_{2}\left(x,y\right)dxdy\right)^{\frac{1}{p'}}\left(\int_{0}^{\infty}\int_{0}^{\infty}K_{3}\left(x,y\right)dxdy\right)^{1-\lambda}\\
 & \leq\left(\int_{0}^{\infty}x^{pA_{1}}F^{p}\left(x;q'\tilde{p},q'\tilde{q}\right)f^{p}\left(x\right)dx\right)^{\frac{1}{q'}}\left(\int_{0}^{\infty}y^{qA_{2}}G^{q}\left(y;p'\tilde{p},p'\tilde{q}\right)g^{q}\left(y\right)dy\right)^{\frac{1}{p'}}\\
 & \ \ \times\left[\left(\int_{0}^{\infty}y^{qA_{2}}G^{q}\left(y;p'\tilde{p},p'\tilde{q}\right)g^{q}\left(y\right)dy\right)\left(\int_{0}^{\infty}x^{pA_{1}}F^{p}\left(x;q'\tilde{p},q'\tilde{q}\right)f^{p}\left(x\right)dx\right)\right]^{1-\lambda}\\
 & \leq\left\{ e^{\frac{p}{q'}\left(q'\tilde{p}+q'\tilde{q}\right)}k^{p}\left(F;q'\tilde{p},q'\tilde{q}\right)\right\} ^{\frac{1}{p}}\left\{ e^{\frac{q}{p'}\left(p'\tilde{p}+p'\tilde{q}\right)}k^{q}\left(G;p'\tilde{p},p'\tilde{q}\right)\right\} ^{\frac{1}{q}}\\
 & \ \ \times\left(\int_{0}^{\infty}x^{\frac{p}{q'}\left(1-s_{1}-s_{2}\right)+p\left(A_{1}-A_{2}\right)}f^{p}\left(x\right)dx\right)^{\frac{1}{p}}\left(\int_{0}^{\infty}y^{\frac{q}{p'}\left(1-s_{1}-s_{2}\right)+q\left(A_{2}-A_{1}\right)}g^{q}\left(y\right)dy\right)^{\frac{1}{q}}.
\end{align*}
where
\[
\left\{ e^{\frac{p}{q'}\left(q'\tilde{p}+q'\tilde{q}\right)}k^{p}\left(F;q'\tilde{p},q'\tilde{q}\right)\right\} ^{\frac{1}{p}}\left\{ e^{\frac{q}{p'}\left(p'\tilde{p}+p'\tilde{q}\right)}k^{q}\left(G;p'\tilde{p},p'\tilde{q}\right)\right\} ^{\frac{1}{q}}=e^{2\left(\tilde{p}+\tilde{q}\right)}k\left(F;q'\tilde{p},q'\tilde{q}\right)k\left(G;p'\tilde{p},p'\tilde{q}\right).
\]

It marks the end of the proving of the first part of the theorem. The equivalence of (105) and (106) and the conditions for equality to hold can be proved in the way which explained in \cite[p. 649, Theorem 1]{Hilbert Inequality and Gaussian Hypergeometric Functions}.
\end{Proof}
\begin{Remark}
In fact, the way we prove the inequalities (105) and (106) provides a common method to generalize inequalities related to gauss hypergeometric functions. Once some integral identities are established, we can prove it according to the primary way of proving. For instance, the inequalities in the reference \cite{Some Extensions of Hilberts Integral Inequality} can be generalize in the same way. (we should notice that we use \textbf{definition 2.2} instead of \textbf{2.1}).

On another aspect,we can directly use the definition of extended hypergeometric functions to construct inequalities.because the new definition is strong and flexible enough to get the desired results.
\end{Remark}

\section{Conclusion}
By using extended Beta functions we have defined some extensions of generalized hypergeometric function and Lauricella's hypergeometric function. In some special cases, those extensions have been studied in \cite{Extended hypergeometric and confluent hypergeometric functions}, \cite{Some generating relations for extended hypergeometric function via generalized fractional derivative operator} and \cite{Extension of gamma beta and hypergeometric functions}. With flexible manipulating of extended definitions, we have obtained some generalization of properties of hypergeometric functions in classical theories. We have also explained an example of Hardy-Hilbert type inequality involving extended Gauss hypergeometric functions in detail, which may be a precedent for applying these extension in other branches of mathematics. It is worthwhile to point out that some fundamental relations between fractional calculus and extended generalized hypergeometric functions have been proposed in this paper. If we want to study it in the aspect of fractional calculus, we need to get more properties of the fractional operator defined by (43). As far as this is concerned, we have gotten some ideas, which will be presented in next publications.

\section*{Acknowledgement}
The author is grateful to the anonymous referee for his/her valuable comments and suggestions on the improvement of this paper.

\end{document}